\documentclass[letter]{amsart}
\usepackage{amsaddr}

\newcommand{\figscale}{.75}

\usepackage{graphicx}
\usepackage{graphics}
\usepackage{amsmath,amssymb}
\usepackage[inline]{enumitem}
\usepackage{xifthen}
\usepackage{tikz-cd}
\usepackage{braket}
\usepackage{mathtools}

\usepackage{amsthm}

\usepackage{mathrsfs}

\usepackage{hyperref}

\setenumerate{leftmargin=*}

\theoremstyle{plain}
\newtheorem{thm}{Theorem}[section]

\newtheorem{lem}[thm]{Lemma}
\newtheorem*{lem*}{Lemma}
\newtheorem{prop}[thm]{Proposition}
\newtheorem{cor}[thm]{Corollary}
\newtheorem*{cor*}{Corollary}

\newtheorem*{claim*}{Claim}
\newtheorem*{thm*}{Theorem}

\newtheoremstyle{break}%
{}{}%
{\itshape}{}%
{\bfseries}{}
{\newline}{}

\theoremstyle{break}

\theoremstyle{definition}
\newtheorem{definition}[thm]{Definition}

\newtheorem*{nota}{Notation}

\theoremstyle{remark}

\newtheorem{obs*}[thm]{Observation}
\newtheorem*{rmk}{Remark}
\newtheorem{rmk*}[thm]{Remark}

\newcommand{\Z}{\mathbb{Z}}
\newcommand{\R}{\mathbb{R}}

\newcommand{\ie}{\textit{i.e.}}
\newcommand{\eg}{\textit{e.g.}}
\newcommand{\be}{\begin{enumerate}}
\newcommand{\ee}{\end{enumerate}}
\newcommand{\bi}{\begin{itemize}}
\newcommand{\ei}{\end{itemize}}
\newcommand{\itembf}[2]{%
\ifthenelse{\isempty{#1}}{\item {\bfseries #2}}{\item[\textbf{#1}]{\bfseries #2}}%
}

\renewcommand{\bar}{\overline}

\newcommand{\dfn}{\mathrel{\mathop:}=}

\DeclareMathOperator{\img}{im}

\DeclareMathOperator{\st}{st}

\DeclareMathOperator{\im}{im}

\DeclareMathAlphabet{\mathbbold}{U}{bbold}{m}{n}
\DeclareMathAlphabet{\cmb}{OT1}{cmbr}{m}{n}

\DeclareMathOperator{\diam}{diam}

\DeclareMathOperator{\End}{End}
\DeclareMathOperator{\supp}{supp}

\DeclareMathOperator{\Isom}{Isom}

\DeclareMathOperator{\SL}{SL}

\DeclareMathOperator{\stab}{stab}

\DeclareMathOperator{\Mod}{Mod}
\DeclareMathOperator{\PMod}{PMod}

\newcommand{\scratch}[1]{{\color{red}\textbf{#1}}}

\newlist{clist}{enumerate}{1}
\setlist[clist]{label=(\roman*),itemsep=0ex,
  topsep=.5ex,parsep=0ex,leftmargin=3.5em}

\newlist{citem}{itemize}{1}
\setlist[citem,1]{label=\textbullet,topsep=1ex,itemsep=1ex,parsep=0ex}

\newcommand{\sqed}{\hfill $/\!\!/$}

\ifdefined\clean
\renewcommand{\scratch}[1]{{}}
\else
\fi

\DeclareMathOperator{\pac}{\mathscr A}
\DeclareMathOperator{\asdim}{asdim}

\begin{document}

\title{Prescribed arc graphs}
\author{Michael C.\ Kopreski}
\email{kopreski@math.utah.edu}

\begin{abstract}
Given a compact surface $\Sigma$ with boundary and a relation
$\Gamma$ on $\pi_0(\partial\Sigma)$, we define the \textit{prescribed
  arc graph} $\pac(\Sigma,\Gamma)$ to be the full subgraph of the arc
graph $\pac(\Sigma)$ containing only classes of arcs between boundary
components in $\Gamma$.  We prove that $\pac(\Sigma,\Gamma)$ is
connected and infinite-diameter (if $\Sigma$ is not the sphere with
three boundary components), and classify when it is Gromov hyperbolic:
in particular, $\pac(\Sigma,\Gamma)$ is Gromov hyperbolic if and only
if $\Gamma$ is not bipartite, except in some sporadic cases.
\end{abstract}

\maketitle

\section{Introduction, definitions and main results}
Given a compact surface $\Sigma$ with boundary, we recall 
the \textit{arc graph} $\pac(\Sigma)$, whose vertices are 
isotopy classes of essential simple arcs and whose edges are 
determined by disjointness up to isotopy.
We propose the \textit{prescribed arc graph} as generalization of the arc graph, defined to be the full subgraph $\pac(\Sigma,\Gamma) \subset \pac(\Sigma)$ spanned by the isotopy classes of arcs between boundary components in $\Gamma$, a symmetric relation on $\pi_0(\partial \Sigma)$:

\begin{definition}\label{def:pac}
  Let $\Sigma$ be a compact, orientable surface with boundary.
  Given $\Gamma$ a graph with $V(\Gamma) = \pi_0 (\partial \Sigma)$,
  let an essential simple arc $\alpha$ be
  \textit{$\Gamma$-allowed} if $\alpha$ terminates on (not necessarily
  distinct) boundary
  components $a_-, a_+$ such that $(a_-,a_+) \in E(\Gamma)$.
  Define the \textit{$\Gamma$-prescribed arc graph} $\mathscr A(\Sigma,
  \Gamma)$ as follows: \[
  \begin{aligned}
    V(\pac(\Sigma,\Gamma)) &= \{ a \text{ an unoriented isotopy class of
      $\Gamma$-allowed arcs} \}
    \\
    E(\pac(\Sigma, \Gamma)) &= \{(a, a') : a, a'
    \text{ have disjoint representatives}\}
  \end{aligned} \]
\end{definition}
\noindent
We will call $\Gamma$ the \textit{prescribing graph} for
$\pac(\Sigma,\Gamma)$.  We note that if $\Gamma$ is the complete graph
with loops on
$\pi_0(\partial \Sigma)$, then $\pac(\Sigma,\Gamma)$ is the usual arc
graph.

Henceforth let $\Sigma_g^b$ denote the
compact, orientable surface with genus $g$ and $b$ boundary
components. We prove some initial results concerning the geometry of $\pac(\Sigma,\Gamma)$, including the following:

\begin{thm}\label{thm:conninf}
  Assume that $\chi(\Sigma) \leq -1, E(\Gamma) \neq \varnothing$, and
  $\Sigma \neq \Sigma_0^3$.  Then $\pac(\Sigma, \Gamma)$ is connected
  and has infinite diameter.
\end{thm}

\begin{thm}\label{thm:hypclass0}
  Assume that $\chi(\Sigma) \leq -1, E(\Gamma) \neq \varnothing$, and
  $\Sigma \neq \Sigma_0^3$.  Then if $\Sigma = \Sigma_0^{n+1}$ and $\Gamma$ is a
$n$-pointed star, or if $\Sigma = \Sigma_1^2$ and $\Gamma$ is a
non-loop edge, then $\pac(\Sigma,\Gamma)$ is $\delta$-hyperbolic. Otherwise, $\pac(\Sigma,\Gamma)$ is $\delta$-hyperbolic if and only if
  $\Gamma$ is not bipartite.
\end{thm}

Section~\ref{sec:conn} proves Theorem~\ref{thm:conninf} outside of
some sporadic low-complexity cases, as well as establishing an upper
bound on distances $d(a,b)$ in $\pac(\Sigma,\Gamma)$
in terms of intersection number $i(a,b)$. Section~\ref{sec:hyp} shows that $\pac(\Sigma,\Gamma)$ is
$\delta$-hyperbolic if $\Gamma$ is not bipartite, again excluding some
sporadic cases.

To complete the proof of Theorem~\ref{thm:hypclass0}, we appeal to the existence of disjoint \textit{witness subsurfaces}.  As with the usual arc graph, we define witnesses to be subsurfaces that
are cut by every $\Gamma$-allowed arc:

\begin{definition}
  An essential connected proper  
subsurface $W \subset \Sigma$ is a
  ($\Gamma$-)\textit{witness} for $\pac(\Sigma, \Gamma)$ if every
  $\Gamma$-allowed arc intersects $W$.
\end{definition}

\noindent
Section~\ref{sec:nonhyp} proves (again,
ignoring some sporadic cases) that $\Gamma$ is bipartite if and only
if 
there exist a pair of disjoint $\Gamma$-witnesses; the latter implies
a quasi-isometric embedding of $\Z^2$.  Finally, Section~\ref{sec:sporadic} addresses the sporadic cases
missing from the preceding sections.  

In the spirit of \cite{diskcpx}, we may expect that
$\pac(\Sigma,\Gamma)$ is $\delta$-hyperbolic if and only if there does
not exist any pair of disjoint $\Gamma$-witnesses.  This hypothesis in
fact holds true, although we do not prove it directly.  From the results in Sections~\ref{sec:hyp} and \ref{sec:nonhyp} along with the sporadic cases in Section~\ref{sec:sporadic}, we conclude:

\begin{thm}\label{thm:hypclass}
  Assume that $\chi(\Sigma) \leq -1, E(\Gamma) \neq \varnothing$, and
  $\Sigma \neq \Sigma_0^3$.  Then
$\pac(\Sigma,\Gamma)$ is
  $\delta$-hyperbolic if and only if $\Sigma$ does not admit two
  distinct, disjoint $\Gamma$-witnesses that are not homeomorphic to
  $\Sigma_0^3$.
\end{thm}

We remark that in the case where $\Gamma$ is loop-free, \ie\ every
$\Gamma$-allowed arc terminates on two distinct boundary components,
the proofs of the above (and in particular, those in Section~\ref{sec:conn})
simplify considerably.  

\subsection{Motivation}
Although perhaps overshadowed in recent decades by the curve 
complex $\mathscr C(\Sigma)$, the arc graph $\pac(\Sigma)$ 
has been an object of intrinsic interest in the classical 
study of surfaces of finite topological type 
(\ie\ compact, with finitely many marked points) 
and their mapping 
class groups.  For example, the arc complex of a finite type surface $\Sigma$ triangulates its Teichm\"uller 
space $\mathcal T(\Sigma)$ \cite{harer}, with natural coordinates arising from horodisks in (cusped) hyperbolic metrics \cite{bowep}; earlier work \cite{mosher} uses $\pac(\Sigma)$ to study conjugacy classes of $\Mod(\Sigma)$.

Nonetheless, the current work arises instead from combinatorial objects defined for surfaces of infinite topological type (\ie\  with
infinite genus or infinitely many
boundary components or marked points, and whose fundamental group is
not finitely generated).
For such surfaces, we have that 
$\diam(\mathscr C(\Sigma)) = 2$, and likewise $\diam(\pac(\Sigma)) = 2$ whenever the number of marked points or boundary components is infinite.
A number of authors have proposed suitable (\eg\ infinite diameter, connected, and
$\delta$-hyperbolic) combinatorial models in the infinite type setting on which $\Mod(\Sigma)$
acts continuously, in analogue to $\mathscr C(\Sigma)$ and $\pac(\Sigma)$ in the
finite type setting.  Danny Calegeri's original 2009 blog post
on 
mapping class groups of infinite type surfaces \cite{danny}
defines the ray graph on $\R^2
\setminus \text{Cantor set}$, which was shown to be infinite diameter,
connected, and $\delta$-hyperbolic by Juliette Bavard in \cite{bavard}.  
More recently, we note the omnipresent arc graph defined by Fanoni, 
Ghaswala, and McLeay \cite{fgm} and the grand arc graph defined 
by Bar-Natan and Verberne
\cite{grand}.  We pay special attention to the grand arc graph 
$\mathscr G(\Sigma)$, in which one
considers only arcs between maximal ends of distinct type, in the
sense of the partial order on $\End(\Sigma)$ defined in \cite{mannrafi}.

One observes that, in the definitions above, the combinatorial
model is made ``sparse enough'' by restricting which arcs in $\Sigma$
are considered.  We propose the prescribed arc graph as similar generalization of the usual arc
graph $\mathscr A(\Sigma)$ in the finite type setting, both inspired and motivated
by the combinatorial models for infinite type surfaces discussed
above.  For example, let $\Sigma$ be a surface of infinite type.  Then
for suitable compact exhaustion by finite type subsurfaces
$\Sigma_i \subset
\Sigma$ and prescribing graphs $\Gamma_i$ on $\pi_0(\partial
\Sigma_i)$,  $\pac(\Sigma_i,\Gamma_i)$ is meant to ``approximate''
$\mathscr G(\Sigma)$.  To start, one may show that for
appropriate $(\Sigma_i,\Gamma_i)$, 
$\pac(\Sigma_i,\Gamma_i)$ coarsely embeds into $\mathscr G(\Sigma)$,
hence \eg\ $\asdim \mathscr G(\Sigma) \geq \asdim
\pac(\Sigma_i,\Gamma_i)$, through which we aim to show that the
asymptotic dimension of $\mathscr G(\Sigma)$ is infinite.  We will discuss
these arguments in a future work.

Finally, we recall the results of Aramayona and Valdez 
\cite{aramayona}, which classify the connectedness, 
diameter, and $\delta$-hyperbolicity for a broad class of 
``sufficiently invariant''  subgraphs of the arc and curve 
graph $\mathscr{AC}(\Sigma)$, for $\Sigma$ a surface of 
infinite type.  Here, arcs are assumed to terminate on a set of marked points $\Pi$, and a subgraph $\mathscr G \subset 
\mathscr{AC}(\Sigma)$ is called \textit{sufficiently 
invariant} if there exists a subset $P \subset \Pi$ such that $\mathscr G$ is preserved setwise by the subgroup of the relative 
mapping class group $\Mod(\Sigma,\Pi)$ that  preserves $P$ setwise.  An extension of our results to arc graphs on infinite type surfaces $\Sigma$ (\eg\ those arising from direct limits of prescribed arc graphs $\pac(\Sigma_i,\Gamma_i)$ on compact surfaces $\Sigma_i$) would likewise extend the work of \cite{aramayona}, in the sense that the resulting arc subgraphs would be invariant only up to subgroups of $\Mod(\Sigma,\Pi)$ that fix a particular symmetric relation $\Gamma$ on $\Pi$.

\subsection{Trivial cases and marked points}  

If $\Sigma$ is the closed annulus or disk, $\pac(\Sigma,\Gamma)$ is
either empty or a singleton.  
We will ignore these trivial cases and assume $\Sigma$ is neither.  If $E(\Gamma) = \varnothing$ or $\Sigma = \Sigma_0^3$, then $\pac(\Sigma,\Gamma)$ is likewise empty or finite; we typically exclude these cases as well.

We note that omission of marked
points or punctures from the definition of the prescribed arc graph is
one of convenience, and likewise the choice that arcs terminate on
boundary components instead of marked points.
In particular, without loss of generality we will
assume every surface considered is without marked points: if $\Sigma$
is a marked surface and $\Sigma'$ is obtained by deleting a neighborhood
of each marked point, then $\pac(\Sigma) \cong \pac(\Sigma')$ and the 
prescribed arc graph may be taken to be the corresponding full subgraph 
of $\pac(\Sigma')$.
Moreover, Definition~\ref{def:pac} could just as well have been taken
with respect to arcs between marked points, with $V(\Gamma)$
corresponding to marked points instead of boundary components.
We will occasionally appeal
to this viewpoint when it is more appropriate.

We also define some convenient notation:

\begin{nota}
Let $\alpha$ be an oriented essential arc between boundary components
of a surface $\Sigma$.  We
denote by $\alpha_-,\alpha_+$ the boundary component containing the
initial, resp.\ terminal endpoint of
$\alpha$.  If $a$ is an isotopy class of such arcs,
then $a_\pm \dfn \alpha_\pm$ for a choice of representative arc $\alpha$.  Let
$\partial \alpha \dfn \{\alpha_\pm\}$ and $\partial a \dfn
\{a_\pm\}$.
\end{nota}

\begin{rmk}
For simplicity, we typically conflate arcs and their
isotopy classes in proofs where well defined.
For example, while $a \in V(\pac(\Sigma,\Gamma))$
denotes an isotopy class of $\Gamma$-allowed arcs, we may apply
topological operations such as intersection or
concatenation, or ascribe properties like transversality and minimal
position,
by (implicitly) choosing a representative $\alpha \in a$.  Conversely,
given a $\Gamma$-allowed arc $\delta$, we may view $\delta \in
V(\pac(\Sigma,\Gamma))$ by (implicitly) passing to the isotopy class $[\delta]$.
\end{rmk}

We conclude this section with an important final definition and some
initial discussion useful for the remaining work.

\subsection{A canonical contraction}

Consider any subgraph $\Gamma' \subset \Gamma$; without loss of
generality, assume $V(\Gamma') = V(\Gamma) = \pi_0(\partial \Sigma)$,
else append singletons for the remaining vertices.  If $a$ is an
isotopy class of $\Gamma'$-allowed arcs, then likewise it is
$\Gamma$-allowed.  Let
$\pi : \pac(\Sigma,\Gamma') \to \pac(\Sigma,\Gamma)$ be the simplicial
map obtained by the inclusion of $V(\pac(\Sigma,\Gamma'))$ into
$V(\pac(\Sigma,\Gamma))$; since $\pi$ is simplicial, it is
$1$-Lipschitz.  Moreover, if $\Gamma$ is loop-free, then every $\Gamma$-allowed arc 
$a \in V(\pac(\Sigma,\Gamma))$ is non-separating, and in particular, 
if $(u,v) \in E(\Gamma')$ then there exists a $\Gamma'$-allowed arc 
between $u$, $v$ disjoint from $a$:

\begin{lem}\label{lem:1onto}
  Suppose that $\Gamma$ is loop-free and $\Gamma' \subset \Gamma$ contains an edge.  Then
  $\pi : \pac(\Sigma,\Gamma') \to \pac(\Sigma,\Gamma)$ is $1$-coarsely
  surjective. \qed
\end{lem}

Finally, we prove a technical lemma of
importance in Sections~\ref{sec:conn} and \ref{sec:hyp}:

\begin{lem}\label{lem:loop_onto}
Suppose that $\ell_0 \subset \Gamma$ is a loop and 
$w \subset \partial \Sigma$ is the boundary component in
$\ell_0$.  Then the vertex set
$X_2 \dfn \{a \in V(\pac(\Sigma,\Gamma) : d(a, \pi\pac(\Sigma,\ell_0))
> 1 \}$ is independent in $\pac(\Sigma,\Gamma)$ and contains only arcs with
an annular complementary component containing $w$.
\end{lem}
\begin{proof}
  It suffices to show that if $a \in X_2$, then $a$
  has an annular complementary component containing $w$: we note that
  any two such non-isotopic arcs must intersect.  Let $\Sigma_\pm$ denote the
  closure(s) of the complementary component(s) of $a$.
Since $a$ is essential,
the canonical maps $\Sigma_\pm \to \Sigma$ are $\pi_1$-injective,
hence any essential arc $\delta \subset \Sigma_\pm$ terminating on
$w \cap \partial \Sigma_\pm$ descends to an essential
$\ell_0$-allowed arc in $\Sigma$ disjoint from $a$, and $d(a,
\pi\pac(\Sigma,\ell_0)) = 1$, a contradiction.  Thus no such $\delta$
exists, hence either $\Sigma_\pm$ are both disks, a contradiction
since $\Sigma$ is not a disk, or one component is an annulus
containing $w$.
\end{proof}
\begin{definition}
We will call the vertex set $X_2$ \textit{exceptional}.
\end{definition}

\section{Connectedness of $\pac(\Sigma,\Gamma)$}
\label{sec:conn}

We first define a generalization of the unicorn arcs introduced in \cite{hpw}:

\begin{definition}
Let $a,b$ be isotopy classes of $\Gamma$-allowed arcs.  Then a
\textit{$\Gamma$-unicorn} of $a,b$ is (the isotopy class of) a
$\Gamma$-allowed concatenation $\alpha_0 * \beta_0$ of subarcs
$\alpha_0 \subset \alpha \in a$ and $\beta_0 \subset \beta \in b$. The
\textit{$\Gamma$-unicorn subgraph} $\mathscr U(a,b) \subset
\pac(\Sigma,\Gamma)$ is the full subgraph spanned by $\Gamma$-unicorns
of $a,b$.\end{definition}

In addition, we impose that $a,b \in \mathscr U(a,b)$; we will say
that a unicorn is \textit{proper} if $x \neq a,b$. Unlike in the usual
arc graph, $\mathscr U(a,b)$ is not always connected. For example,
since any proper unicorn between $a, b$ has one endpoint in $\partial
a$ and the other in $\partial b$, if $\partial a$ shares no edges with
$\partial b$ in $\Gamma$ then no such arc is $\Gamma$-allowed: if
$a,b$ intersect, then $\mathscr U(a,b)$ consists of two singletons.
In fact this is the only case for which $\mathscr U(a,b)$ is
disconnected:

\begin{prop}\label{prop:unicorn}
Let $a, b \in V(\pac(\Sigma,\Gamma))$.  If $\partial a$ and $\partial
b$ share an edge in $\Gamma$, then $\mathscr U(a,b)$ is connected.
Moreover, $d(a,b) \leq i(a,b) + 1$, where $i$ denotes the geometric
intersection number.
\end{prop}

\begin{proof}
  Let $x \in \mathscr U(a,b)$; up to exchanging $a,b$, assume $x \neq
  a$.  Orient $a, b$ and $x$ such that $b_+ = x_+$ and $(a_-,b_+) \in
  E(\Gamma)$.  Assume $a,b$ are in minimal position;
  if $x = b$, then let $\beta_0 = b$, else fix subarcs $\alpha_0 \subset a$ and $\beta_0 \subset b$ such that $x = \alpha_0 * \beta_0$.  It suffices to find a $\Gamma$-unicorn $\gamma$ of $a,b$ disjoint from $x$ with $i(\gamma,a) <
  i(x,a)$; by induction on intersection number, we conclude.


  Let $s \in a \cap x$ be the first transverse intersection with $x$
along $a$; note $s \in \beta_0$.  Let $\gamma$ be the
  concatenation of the subarc from $a_-$ to $s$ along $a$,
  and the subarc from $s$ to $x_+$ along $x$; the latter subarc lies
  in $\beta_0$ and $\partial \gamma = (a_-,b_+) \in E(\Gamma)$, hence
  if $\gamma$ is essential, then it is a $\Gamma$-unicorn.  However,
  if $\gamma$ bounds a half-disk, then $a,
  b$ share a bigon or half-bigon, hence are not in
  minimal position. Finally, we observe that (up to isotopy) $\gamma,
  x$ are disjoint and that $i(\gamma,a) \leq i(x,a) - 1$.
\end{proof}

Despite the non-existence of some unicorn paths, we may now show that $\pac(\Sigma,\Gamma)$ is connected if $\Gamma$ is loop-free.  In particular, it suffices that $\pac(\Sigma, e)$ is connected for some
edge $e \in \Gamma$: $\pi
\pac(\Sigma,e) \subset \pac(\Sigma,\Gamma)$ is thus a connected,
$1$-dense subgraph by Lemma~\ref{lem:1onto}. Observing that for
any edge $e$ and $a,b \in V(\pac(\Sigma,e))$ $\partial a = \partial b$, $\pac(\Sigma,e)$ is connected by Proposition~\ref{prop:unicorn}.

If $\Gamma$ contains a loop $\ell_0$, a similar argument
\textit{almost} applies: by Proposition~\ref{prop:unicorn}
$\pac(\Sigma,\ell_0)$ is connected, and by Lemma~\ref{lem:loop_onto}
$V(\pac(\Sigma,\Gamma)) \setminus X_2$ lies within the connected
$1$-neighborhood of $\pi(\Sigma,\ell_0)$.  However, to show that the
exceptional vertices $X_2$ are not isolated we must defer to
Section~\ref{sec:d_bound}, where we prove the distance estimate in Proposition~\ref{prop:unicorn} in general.  This result will be useful in
addition for our discussion in Section~\ref{sec:hyp}.

Finally, we observe that the proof of Proposition~\ref{prop:unicorn}
implies the following:
\begin{cor}\label{cor:or_uni}
  Suppose $\alpha, \beta$ are oriented $\Gamma$-allowed arcs such that
  $\alpha_+$ and $\beta_+$ lie in the same boundary component of
  $\Sigma$.  If $\alpha, \beta$ are in minimal position, then any
  simple arc obtained by oriented surgery of $\alpha, \beta$ is a $\Gamma$-unicorn. \qed
\end{cor}

\subsection{An upper bound for distance}\label{sec:d_bound}

We will prove that except in certain low-complexity cases, $d(a,b)
\leq i(a,b) + 1$, extending the result of
Proposition~\ref{prop:unicorn} for generic arcs.
 We
begin with arcs with one or two intersections.

\begin{lem}\label{lem:i_1}
  Assume that $\Sigma \neq \Sigma_0^3$ and that if
  $\Sigma = \Sigma_1^2$, then $\Gamma$ contains a non-loop edge.  If
  $a, b$ are unoriented isotopy classes of $\Gamma$-allowed arcs and
  $i(a,b) = 1$, then $d(a,b) = 2$ in $\pac(\Sigma,\Gamma)$.
\end{lem}

\begin{proof}
  Assume $a,b$ are in minimal position.  If $\partial a \cap \partial b \neq
  \varnothing$ then $\partial a, \partial b$ share an edge in $\Gamma$, and likewise by assumption if $\Sigma = \Sigma_1^2$: by Proposition~\ref{prop:unicorn}, $d(a,b) \leq
  i(a,b) + 1 = 2$.  
  We suppose neither holds and 
consider three cases:
\begin{enumerate}[label=\textit{(\roman*)}]
\item \textit{Neither $a$ nor $b$ are loops (\ie\ $a_+ \neq a_-, b_+ \neq
  b_-$).}  Let $N$ be a regular neighborhood of $a \cup b \cup
  b_\pm$, and let $\delta$ be a component of $\partial
  N$.\footnote{Remark: $\partial N$ denotes the
  \textit{topological} boundary, hence does not contain $b_\pm$ or any
  subarcs of $a_\pm$.}  Then $\delta$ is essential and
  disjoint from $a, b$; since $\delta_\pm = a_\pm$,
  $\delta$ is $\Gamma$-allowed.

  \item \textit{(Without loss of generality) $a$ is a loop but
      $b$ is not ($a_+ = a_-, b_+ \neq b_-$).} We proceed as above, and let $\delta_1,
    \delta_2$ denote the components $\partial N$. 
If either is essential, then we conclude as
above; else both bound half-disks and $\Sigma = \Sigma_0^3$.

    \item \textit{Both $a$ and $b$ are loops ($a_+ = a_-,
        b_+=b_-$).} Let $N$ be a regular neighborhood of $a \cup
      b \cup b_+$.  Since
      $a, b$ intersect once, $N \cup N(a_+)$ is a
      thrice-punctured torus; let $\delta = \partial N$. 
      If $\delta$ is essential, then we conclude
      as above; else $\delta$ bounds a half-disk and $\Sigma =
      \Sigma_1^2$. 
    \end{enumerate}
It follows that $d(a,b) = 2$.
\end{proof}

\begin{lem}\label{lem:i_2}
Assume that $\Sigma \neq \Sigma_0^3$ and that if $\Sigma =
  \Sigma_1^2$, then $\Gamma$ contains a non-loop edge.  If $a, b$ are unoriented isotopy classes
of $\Gamma$-allowed arcs such that $i(a,b) = 2$,
then $d(a,b) \leq 3$ in $\pac(\Sigma,\Gamma)$.
\end{lem}
\begin{proof}
Assume $a,b$ are in  
minimal position and that $\partial a \cap \partial b =
  \varnothing$, else conclude by Proposition~\ref{prop:unicorn} as
  above.  Let $\{s,t\} = a \cap b$ and let 
  $\alpha_0, \beta_0$ denote the subarcs of $a, b$
  respectively between $s$ and $t$; orient $a$ such that 
  $s$ is nearer $a_-$, and likewise with $b$.  Let $\alpha^\pm \subset
  a \setminus \mathring \alpha_0$ denote the subarc
  incident on $a_\pm$, and likewise for $\beta^\pm 
\subset b \setminus \mathring \beta_0$.
 We consider two cases:
\begin{enumerate}[label=\textit{(\roman*)}]
\item \textit{Not both $a, b$ are loops, or $\hat
    \imath(a,b) = \pm 2$, where $\hat \imath$ denotes the
    algebraic intersection number.}  Without loss of generality,
  assume $b$ is a loop only if $a$
  is a loop.  Let $\delta$ be the concatenation of $\beta^-,
  \alpha_0$, and $\beta^+$.  We observe that $\delta$ and $
b$ are disjoint up to isotopy and $i(\delta,a) 
\leq 1$, with equality when $\hat \imath (a,b) = \pm 2$.  Hence if
$\delta$ is essential the claim is shown: since $\partial 
\delta = \partial b$, $\delta$
is $\Gamma$-allowed and $d(\delta, a) \leq 2$ by 
Lemma~\ref{lem:i_1}, hence $d(a,b) \leq d(a,\delta) + d(\delta, 
b) \leq 3$.  Finally, to verify that $\delta$ is essential, we observe
that either $b$ (hence $\delta$) is not a loop or both $a,b$ are loops
and $\hat\imath(a,b) = \pm 2$, hence $i(\delta,a) = 1$.  In both
cases, it follows that $\delta$ is non-separating.
  \item \textit{Both $a, b$ are loops and $\hat \imath(a,b) = 0$.}
    Let $N_a$ be a regular neighborhood of $a \cup b
    \cup b_\pm$ and $N_b$ a regular neighborhood of $a
    \cup b \cup a_\pm$.  Let $\mathscr D_a$ be the
    set of arc components of $\partial N_a$, and
    likewise $\mathscr D_b$ for $\partial N_b$.  For any $\rho \in \mathscr D \dfn \mathscr D_a \cup \mathscr
    D_b$, $\rho$ is disjoint from $a$, $b$ and either
    $\partial \rho \in \partial a$ or $\partial \rho \in \partial b$.
    Hence if $\rho$ is essential, then $\rho$ is $\Gamma$-allowed and
    $d(a,b) = 2$.

    Assume every arc in $\mathscr D$ bounds a half-disk.  Then (\eg\ by an Euler characteristic argument) the
    simple closed curve $\alpha_0 \cup \beta_0$ separates $\Sigma =
    \Sigma_0^3 \cup_{\alpha_0 \cup \beta_0} \Sigma'$,
where
    $\Sigma'$ is the subsurface with boundary
    $\alpha_0 \cup \beta_0$ not containing $a_\pm$ and $b_\pm$.
    If $\Sigma'$ is a disk, then $a, b$ are not in minimal
    position, and if $\Sigma'$ is an annulus, then $\Sigma = \Sigma_0^3$,
    both contradictions with our assumptions.  Hence there
    exists a simple arc $\gamma \subset \Sigma'$ between $s, t \in \partial
    \Sigma'$ that does not
    bound a half-disk in $\Sigma'$.  Let $\alpha' = \alpha^- * \gamma *
    \alpha^+$ and $\beta' = \beta^- * \gamma * 
\beta^+$.  If $\alpha'$ bounds a half-disk in $\Sigma$, then
    likewise does $\gamma$ in $\Sigma'$, and similarly with $\beta'$:
    both $\alpha' $ and $\beta'$ are essential.  Since $\partial
    \alpha' = \partial a$ and $\partial \beta' = \partial b$,
    $\alpha', \beta'$ are $\Gamma$-allowed, and we observe that
    $a$ and $\alpha'$, $\alpha'$ and $\beta'$, and $\beta'$ and
    $b$ are disjoint up to isotopy.  $d(a,b) \leq 3$. \qedhere
\end{enumerate}
\end{proof}

\begin{figure}[t]
\begin{tikzpicture}
\pgftext{%
\includegraphics[scale=\figscale]{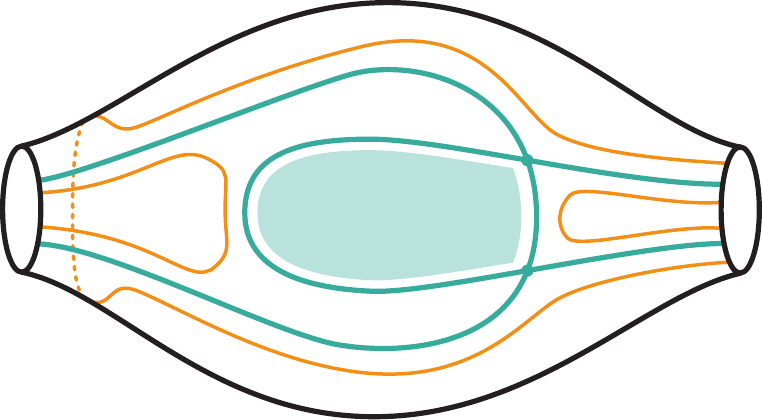}%
};
\node at (0,0) {\Large $\Sigma'$};
\node at (3.3,0) {$b_\pm$};
\node at (-3.3,0) {$a_\pm$};

\node at (1.4,.6) {$s$};
\node at (1.4,-.7) {$t$};

\end{tikzpicture}
\caption{A maximal collection of disjoint arcs (orange, thin) in $\mathscr D$, each bounding a half disk.}
\end{figure}
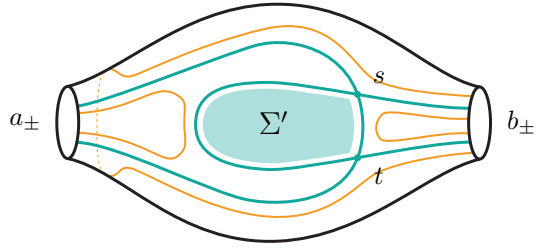

For $i(a,b) > 2$, we construct a sequence of arcs with strictly
decreasing intersection number and proceed by induction. 

\begin{prop}\label{prop:dist_i}
  Assume that $\Sigma \neq \Sigma_0^3$ and that if $\Sigma =
  \Sigma_1^2$, then $\Gamma$ contains a non-loop edge.   If $a, b$ are
  unoriented isotopy classes of $\Gamma$-allowed arcs and $i(a,b) = n$, then
  $d(a,b) \leq n + 1 $
  in $\pac(\Sigma,\Gamma)$. 
\end{prop}
\begin{proof}
  As above, we may assume that $\partial a \cap \partial b =
  \varnothing$, else conclude by Proposition~\ref{prop:unicorn}.
  Assume $a,b$ are in minimal position.  It
  suffices to find a $\Gamma$-allowed arc $\delta$
  disjoint from $b$ such that  $i(\delta , a) < i(a,
  b)$, or likewise exchanging $a$ and $b$; 
  we then conclude by
  induction on intersection number, noting that the cases $i(a,b) \leq
  2$ follow from Lemmas~\ref{lem:i_1} and \ref{lem:i_2}.  Assume
  $i(a,b) \geq 3$ and let $s
  \in a \cap b$ be the first intersection
  along $b$; let $t$ be the subsequent intersection along
  $a$.   Let 
  $\alpha_0, \beta_0$ denote the subarcs of $a, b$
  respectively between $s$ and $t$.  Let $\beta^\pm \subset
  b \setminus \mathring \beta_0$ denote the subarc
  incident on $b_\pm$. We consider two cases:
  \begin{enumerate}[label=\textit{(\roman*)}]
  \item \textit{$s,t$ are intersections of the same sign.}  Let
    $\delta = \beta^- * \alpha_0 * \beta^+$.  Up to
    isotopy $\delta$ and $b$ are disjoint and $i(\delta,a)
    \leq i(a, b) - 1$.  Since $ \delta_\pm= 
    b_\pm$, if $\delta$ is essential then it is $\Gamma$-allowed.
    If $b$ is not a loop, then $\delta$ is not a loop and hence
    essential.  Suppose that $b$ is a loop.  Isotoping $\delta$ to be
    disjoint from $b$, we observe that the endpoints of $\delta$
    separate the endpoints of $b$ on $b_\pm$:
    $\delta$ is non-separating and 
    hence essential.

\item \textit{$s,t$ are intersections of opposite sign.}  Let
    $\delta$ as above; $\delta$
    and $b$ are disjoint up to isotopy and $i(\delta,a)
    \leq i(a, b) - 2$, hence if $\delta$ is essential,
    then we conclude.  Suppose instead that $\delta$ bounds a half-disk $D$, hence
    $\delta$ and $b$ are loops. Then since $\beta^- \cap
    a = \{s\}$, we have $\beta^+ \cap a =
    \{t\}$, else $a,b$ share a bigon and are not in minimal
    positon. Let $N$ be a
    regular neighborhood of $a \cup D \cup b_\pm$, and let $\xi$
    be the arc component of $\partial N$ not parallel to $a$.
    We observe that $\xi$ is
    disjoint from $a$ and that $i(\xi,b) \leq i(a,b)
    - 2$;  
    $\partial \xi = \partial a$, hence if $\xi$ is
    essential, then $\xi$ is $\Gamma$-allowed.  
    Suppose that $\xi$
    bounds a half-disk, and observe that $b$ intersects $\xi$ only
    on subarcs parallel to $a$.  Hence $b$ and $\xi$ are
    disjoint, else $b$ shares a bigon or half bigon with $\xi$ incident only on subarcs in $\xi \cap a$, thus also with $a$. 
    Since $a \setminus \mathring \alpha_0$ is parallel to $\xi$,
    $a$ intersects $b$ only at $s,t$, a contradiction since
    $i(a,b) \geq 3$. \qedhere
  \end{enumerate}
\end{proof}

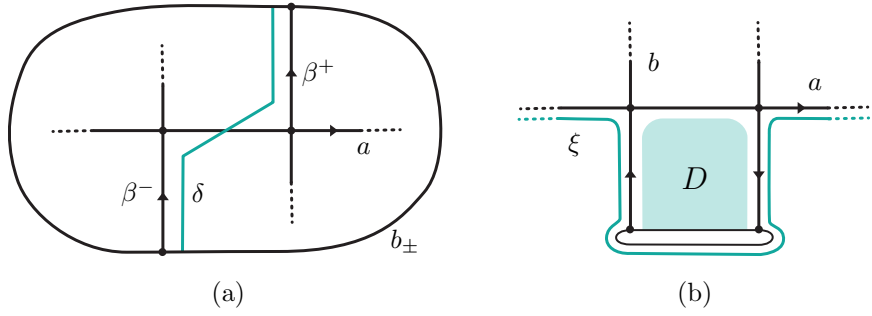
\begin{figure}[t]
\begin{tikzpicture}
\pgftext{%
\includegraphics[scale=.75]{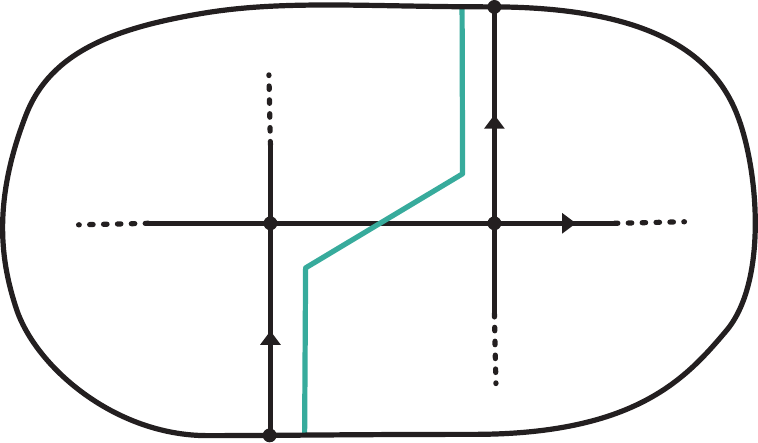}%
\hspace{1cm}%
\includegraphics[scale=.75]{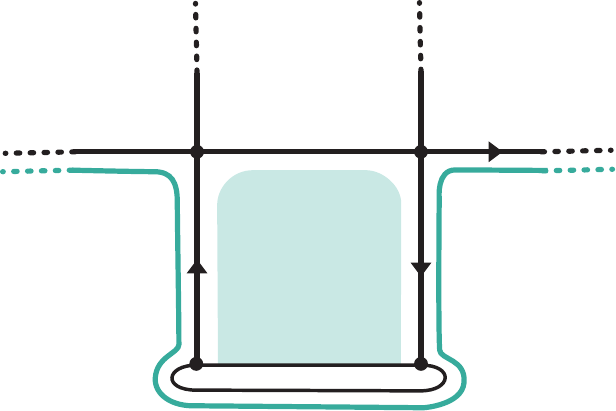}%
}%
\node at (-1,-.25) {$a$};
\node at (-1.6,.75) {$\beta^+$};
\node at (-4,-.85) {$\beta^-$};
\node at (-3.2,-.85) {$\delta$};

\node at (-.45,-1.5) {$b_\pm$};

\node at (-2.8,-2.2) {(a)};
\node at (3.4,-2.2) {(b)};

\node at (3.4,-.65) {\Large $D$};
\node at (5,.6) {$a$};

\node at (2.85,.9) {$b$};
\node at (1.8,-.2) {$\xi$};
\end{tikzpicture}
\caption{(a) $\delta$ separates the endpoints of $b$ in case \textit{(i)}.  (b) the arc $\xi$ in case \textit{(ii)} when $\delta$ bounds a half disk. }
\end{figure}

\begin{rmk*}
Note that if $\Gamma$ is loop-free, then in the statement of the proposition and the two preceding lemmas we may omit our assumptions on $(\Sigma,\Gamma)$. 
\end{rmk*}

Connectivity of $\pac(\Sigma,\Gamma)$ follows as an immediate corollary:
\begin{thm}\label{thm:connected}
  Assume that either $\Gamma$ is loop-free, or
  $\Sigma \neq \Sigma_0^3$ and if $\Sigma = \Sigma_1^2$ then $\Gamma$
  contains a non-loop edge.  $\pac(\Sigma,\Gamma)$ is connected. \qed
\end{thm}
\begin{rmk*}
  The connectedness of $\pac(\Sigma,\Gamma)$ allows us to elaborate
  Lemma~\ref{lem:loop_onto}.  Suppose $\ell_0 \subset \Gamma$ is a
  loop and $X_2\subset V(\pac(\Sigma,\Gamma))$ denotes the set of
  exceptional vertices not in the $1$-neighborhood of
  $\pi\pac(\Sigma,\ell_0)$, and additionally assume
  $\Sigma \neq \Sigma_0^3$ and if $\Sigma = \Sigma_1^2$ then $\Gamma$
  contains a non-loop edge. We observe that for any $v \in X_2$,
  $d(v, \pi\pac(\Sigma,\ell_0)) = 2$.  In particular, since $X_2$ is
  discrete in $\pac(\Sigma,\Gamma)$ by Lemma~\ref{lem:loop_onto} and
  $\pac(\Sigma,\Gamma)$ is connected, $v$ must be adjacent to
  $a \notin X_2$, which lies in the $1$-neighborhood of
  $\pi \pac(\Sigma,\ell_0)$ by the definition of $X_2$. \sqed
\end{rmk*} 

Proposition~\ref{prop:dist_i} also allows us to generalize Lemma~\ref{lem:1onto}:

\begin{lem}\label{lem:3onto}
  Suppose that $\Gamma' \subset \Gamma$ contains an edge, and let
  $\pi : \pac(\Sigma,\Gamma') \to \pac(\Sigma,\Gamma)$ denote the
  simplicial map induced by the inclusion
  $\Gamma' \hookrightarrow \Gamma$.  Then:
\begin{enumerate}[label=(\roman*)]
\item\label{item:3onto:1} if
  $\Gamma$ is loop-free, then $\pi$ is $1$-coarsely surjective.
\end{enumerate}
Moreover, if $\Sigma \neq \Sigma_0^3$ and $\Sigma =
  \Sigma_1^2$ only if $\Gamma$ contains a non-loop edge, then:
\begin{enumerate}[resume, label=(\roman*)]
\item\label{item:3onto:2} if $\Gamma'$ has a non-loop edge, then $\pi$
  is $2$-coarsely surjective; and otherwise
\item\label{item:3onto:3} $\pi$ is $3$-coarsely
  surjective.
\end{enumerate}
\end{lem}
\begin{proof}
  \ref{item:3onto:1} is Lemma~\ref{lem:1onto}.  By change of coordinates, for any
  $a \in V(\pac(\Sigma,\Gamma)) \setminus V(\pi \pac(\Sigma,\Gamma'))$, 
  there exists a $\Gamma'$-allowed arc $\rho$ that
  intersects $a$ at most once if $\rho$ is not a loop, and at
  most twice if $\rho$ is a loop.  We conclude by
  Proposition~\ref{prop:dist_i}.
\end{proof}


\begin{thm}\label{thm:diam}
 Assume $\Sigma \neq \Sigma_0^3$.
  If $\pac(\Sigma)$ has infinite diameter, then likewise does 
  $\pac(\Sigma,\Gamma)$. \qed
\end{thm}

\section{$(\Sigma, \Gamma)$ for which $\pac(\Sigma,\Gamma)$ is hyperbolic}\label{sec:hyp}

In the following, we assume
 that $\Sigma \neq \Sigma_0^3$ and that if 
$\Sigma =  \Sigma_1^2$, 
then $\Gamma$ contains a non-loop edge.  To show the 
hyperbolicity of $\pac(\Sigma,\Gamma)$, 
we apply the following proposition of \cite{bowditch}:
\begin{thm}[Guessing geodesics lemma]\label{thm:gg}
Suppose that $\Omega$ is a connected simplicial graph and there exists
$\{\mathscr U_{a,b}\}_{a,b \in V(\Omega)}$ a family of connected
subgraphs such that $a, b \in \mathscr U_{a,b}$ and $\Delta \geq 0$
such that
\begin{enumerate}[label=(\Roman*)]
\item\label{item:gg1} for $a,b \in V(\Omega)$, if $d(a,b) \leq 1$,
  then $\diam \mathscr U_{a,b} \leq \Delta$, and
\item\label{item:gg2} for $a, b, c\in V(\Omega)$,
  $\mathscr U_{a,c} \subset N_{\Delta}(\mathscr U_{a,b} \cup \mathscr
  U_{b,c})$.
\end{enumerate} 
Then $\Omega$ is $\delta$-hyperbolic for some
$\delta = \delta(\Delta) \geq 0$. \qed
\end{thm}
\noindent
By Proposition~\ref{prop:unicorn}, for $a, b \in V(\pac(\Sigma, \Gamma))$, the $\Gamma$-unicorn subgraph $\mathscr U(a,b)$ is connected only when
$\partial a $ and $\partial b$ share and edge in $\Gamma$.
 Our subgraphs will be
chosen instead to be \textit{augmented unicorn subgraphs} $\mathscr U^+(a, b)$
which lie within a uniformly bounded distance from some connected $\Gamma$-unicorn
subgraph $\mathscr U(a',b')$.

We proceed by first proving $\delta$-hyperbolicity in the case when
$\Gamma$ contains a loop, and then the loop-free case when $\Gamma$
contains an odd cycle.
\begin{prop}\label{prop:loop_hyp}
  Assume that $\Sigma \neq \Sigma_0^3$ and that if
  $\Sigma = \Sigma_1^2$, then $\Gamma$ contains a non-loop edge.  If
  $\Gamma$ contains a loop, then $\pac(\Sigma,\Gamma)$ is uniformly
  $\delta$-hyperbolic.
\end{prop}
\begin{prop}\label{prop:odd_hyp}
  If $\Gamma$ is loop-free and contains an odd cycle, then
  $\pac(\Sigma,\Gamma)$ is uniformly $\delta$-hyperbolic.
\end{prop}
\noindent
Observing that $\Gamma$ is bipartite if and only if it is loop-free
and contains no odd cycles, we conclude the following:
\begin{thm}\label{thm:hyp}
  Assume that either $\Gamma$ is loop-free, or
  $\Sigma \neq \Sigma_0^3$ and if $\Sigma = \Sigma_1^2$ then $\Gamma$
  contains a non-loop edge.  If $\Gamma$ is not bipartite, then
  $\pac(\Sigma, \Gamma)$ is uniformly $\delta$-hyperbolic. \qed
\end{thm}

\subsection{$\Gamma$ contains a loop}

We first consider the case when $\Gamma$ is a single loop, in which
case all $\Gamma$-allowed arcs share a boundary component and the
usual unicorn subgraphs suffice.

\begin{lem}\label{lem:loop_slim}
  If $\Gamma$ is
  comprised of a single, looped edge, then $\{\mathscr U(a,b)\}_{a,b}$ form
  $1$-slim triangles in $\pac(\Sigma, \Gamma)$.
\end{lem}
\begin{proof}
  Let $x \in \mathscr U(a,b)$, and  suppose $x$ is
  comprised of subarcs $\alpha' \subset a, \beta' \subset b$
  with $x_- = \alpha'_- = a_-$ and
  $x_+ = \beta'_+ = b_+$.  Let
  $c \in V(\pac(\Sigma,\Gamma))$, and assume $c$ is in
  minimal position with $x$.

  It suffices to find
  $y \in \mathscr U(a,c) \cup \mathscr U(b,c)$ disjoint
  from $x$.  Without loss of generality suppose that $c$ last
  intersects $x$ at $s \in \alpha'$.  Let $y$ be the
  concatenation of the subarc from $x_-$ to $s$ along $a$ and
  the subarc from $s$ to $c_+$ along $c$. Since $c$ and
  $x$ are in minimal position and $x_+, c_+ = b_\pm$,
  by Corollary~\ref{cor:or_uni} $y \in \mathscr U(a,c)$ is
  $\Gamma$-allowed.  $y$ is disjoint with $x$ up to isotopy.
\end{proof}

We note that if $d(a,b) = 1$ then $a, b$ are disjoint, hence
$\mathscr U(a,b) = \{a,b\}$ has diameter $1$ and Theorem~\ref{thm:gg}
is satisfied for $\Delta = 1$.

\begin{cor}\label{cor:loophyp}
  If $\Gamma$ is comprised of a single, looped edge, then
  $\pac(\Sigma,\Gamma)$ is $\delta$-hyperbolic. \qed
\end{cor}

In fact, it suffices that $\Gamma$ contain a loop.  We assume that
$\Sigma \neq \Sigma_0^3$ and if $\Sigma = \Sigma_1^2$, then $\Gamma$
contains a non-loop edge.  Let $\ell_0 \subset \Gamma$ be a loop in
$\Gamma$ and recall that by Lemma~\ref{lem:loop_onto} the set of
exceptional vertices $X_2$ are discrete in $\pac(\Sigma,\Gamma)$.  Let
$\mathscr A'$ be the metric graph obtained from $\pac(\Sigma,\Gamma)
\setminus X_2$, with the usual graph metric, by adjoining edges of
length $2$ between all neighbors of a vertex $v$ for each $v \in X_2$.
Since $\mathscr A'$ and $\pac(\Sigma,\Gamma)$ are isometric outside of
a collection of disjoint uniformly bounded subsets, they are
uniformly quasi-isometric.  We note that $\pi\pac(\Sigma,\ell_0)$ embeds
isometrically as a $1$-dense subgraph
in
$\mathscr A'$.

\begin{figure}[h]
\begin{tikzpicture}
\pgftext{%
\includegraphics[scale=\figscale]{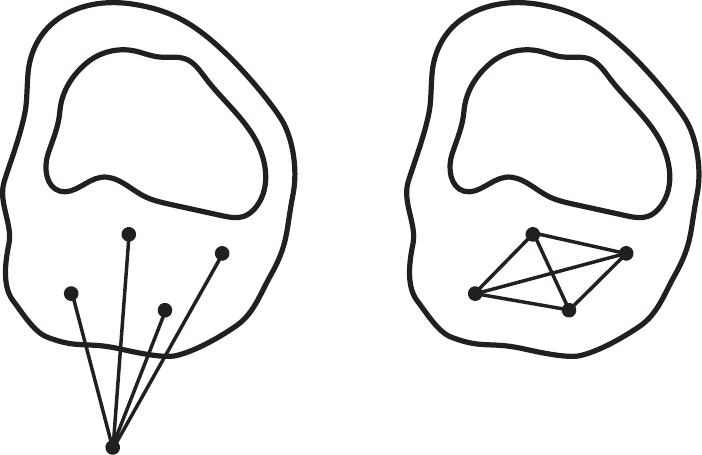}%
}%

\node at (1.6,.8) {$\im \pi$};
\node at (-1.5,.8) {$\im \pi$};

\node at (-1.1,-1.7) {$v \in X_2$};
\node at (-2.1,-.2) {\tiny $N(v)$};

\node at (.95,-.15) {\small $2$};


\node at (1.6,-1.6) {\Large $\mathscr A'$};

\end{tikzpicture}
\caption{Obtaining $\mathscr A'$ from $\pac(\Sigma,\Gamma)$.}
\end{figure}

We prove that $\mathscr A'$ is $\delta$-hyperbolic.  Let
$\xi : \mathscr A' \to \pi\pac(\Sigma,\ell_0)$ be a choice of nearest
point projection.  For $a, b \in V(\pac(\Sigma,\Gamma))$, let
 \[
\mathscr U^+(a,b) \dfn [a,\xi a] \cup \mathscr U_{\ell_0}(\xi a, \xi b) \cup [\xi b, b]
\]
where \eg\ $[a,\xi a]$ denotes an edge and $\mathscr U_{\ell_0}$
denotes the $\pi$-image in $\mathscr A'$ of the respective $\ell_0$-unicorn subgraph.  Since $\pi$ is a contraction,
Lemma~\ref{lem:loop_slim} implies that triangles in
$\{\mathscr U_{\ell_0}(\xi a, \xi b)\}_{a,b}$ are likewise $1$-slim,
and thus for the augmented unicorns $\mathscr U^+(a,b)$ as well.  To
apply Theorem~\ref{thm:gg}, we need only check that \ref{item:gg1} is
satisfied:
\begin{lem}\label{lem:loopcohere}
For any disjoint $a,b \in V(\mathscr A')$, $\diam \mathscr U^+(a,b) \leq 10$.
\end{lem}
\begin{proof}
Let $w$ denote the boundary component in $\ell_0$ and $a' = \xi a, b'
= \xi b$. Since $d(a,a'),d(b,b') \leq 1$, $a,a'$ and $b,b'$ are disjoint up to isotopy.  If both
$a,b \in \pi \pac(\Sigma,\ell_0)$, then $a' = a$ and $b' = b$ are disjoint:
$\mathscr U_{\ell_0}(a',b') = \{a,b\}$.  Assume not, and without loss
of generality, let $b \notin \pi \pac(\Sigma,\ell_0)$.  Then $\partial b \neq \{w\}$; 
orient $b$ such that $b_- \neq w$.  It suffices to show that for any 
$x \in \mathscr U_{\ell_0}(a',b')$, 
$d(x,a) \leq 5$ in $\mathscr A'$.  

Let $s \in b \cap
x$ be the first intersection along $b$; since $b, b'$ are disjoint, $s
\subset a'$.  Let $b_0$ denote the subarc of $b$ between $b_-$ and
$s$, and let $a_0'$ denote the subarc of $x$ along $a'$ between
$a_\pm'$ and $s$.
Let $N$ be a regular neighborhood of $a_0' \cup b_0 \cup
b_-$.  Then $\delta = \partial N$ intersects $x$ at most once and must
be essential, else $\Sigma$ is an annulus:
$d(x,\delta) \leq 2$.  Since $\partial \delta = \{w\}$, $\delta$ is
$\Gamma$-allowed and by Lemma~\ref{lem:loop_onto}
does not lie in $X_2$: $\delta \in V(\mathscr
A')$.  Finally, since $a'$ is disjoint from $a$, $\delta$ intersects
$a$ only if $a$ terminates on $b_-$, hence at most twice.  Applying
Proposition~\ref{prop:dist_i}, $d(a, \delta)
\leq 3$, hence $d(x,a) \leq 5$ as required.
\end{proof}

\begin{figure}[t]
\begin{tikzpicture}
\pgftext{%
\includegraphics[scale=\figscale]{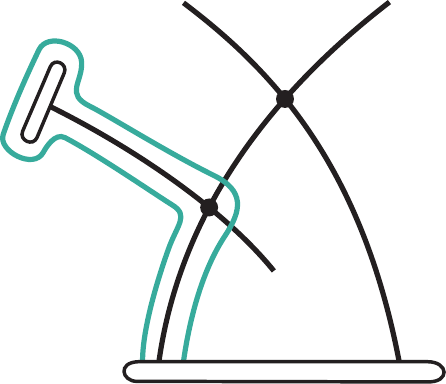}%
}%

\node at (.6,-.6) {$b$};
\node at (1.4,-.2) {$b'$};
\node at (1.3,1.1) {$a'$};

\node at (-.9,-.05) {$\delta$};

\node at (2,-1.35) {$w$};

\end{tikzpicture}
\caption{The arc $\delta$ in the proof of Lemma~\ref{lem:loopcohere}.}
\end{figure}
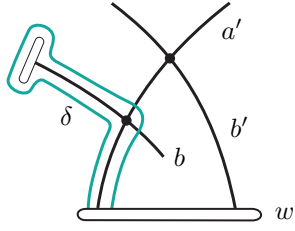

The proof of Proposition~\ref{prop:loop_hyp} follows from
Theorem~\ref{thm:gg} and that
$\pac(\Sigma,\Gamma) \simeq_{\text{q.i.}} \mathscr A'$.

\subsection{$\Gamma$ contains an odd cycle}

Assume that $\Gamma$ is loop-free and contains an odd cycle $C_n$ of
length $n = 2k + 1$.  Fix $e_0 \in E(C_n)$.  By Lemma~\ref{lem:3onto},
the canonical map $\pi : \pac(\Sigma, e_0) \to \pac(\Sigma,\Gamma)$ is
$1$-coarsely surjective.  As above, let
$\xi : \pac(\Sigma,\Gamma) \to \pi \pac(\Sigma, e_0)$ be a choice of
nearest point projection, and for $a,b \in V(\pac(\Sigma,\Gamma))$,
let
\[\mathscr U^+(a,b) \dfn [a,\xi a] \cup \mathscr U_{e_0}(\xi a, \xi b)
  \cup [\xi b, b].\] We verify first that triangles in
$\{\mathscr U_{e_0}(a,b)\}_{a,b \in V(\pi\pac(\Sigma,e_0))}$ are slim
in $\pac(\Sigma,C_n)$; since $\pi$ factors through the $1$-Lipschitz
induced map $\pac(\Sigma,C_n) \to \pac(\Sigma,\Gamma)$, the augmented
unicorns $\mathscr U^+(a,b)$ likewise form slim triangles in
$\pac(\Sigma,\Gamma)$, satisfying \ref{item:gg2} of
Theorem~\ref{thm:gg}.

\begin{lem}\label{lem:hyp_thin}
  Triangles in
  $\{\mathscr U_{e_0}(a,b)\}_{a,b \in V(\pi\pac(\Sigma,e_0))}$ are
  $4$-slim in $\pac(\Sigma,C_n)$.
\end{lem}
\begin{proof}
Let $e_0 = (w_1,w_n)$, and let $w_i$ denote the remaining
boundary components in $C_n$, ordered by adjacency.  Let $a,b,c \in
V(\pi\pac(\Sigma,e_0))$, and let $x \in \mathscr U_{e_0}(a,b)$.  It
suffices to find $y \in \mathscr U_{e_0}(a,c) \cup \mathscr
U_{e_0}(c,b)$ such that $d(x,y) \leq 4$.  Orient $a,b,c$ such that their
initial (resp.\ terminal) points lie in $w_1$ (resp.\ $w_{n}$).
If $x = a,b$, then we may conclude with $y = a, b$
respectively.
Hence without loss of generality, let $x$ be represented by the
concatenation of subarcs $\alpha_0 \subset a$ and $\beta_0
\subset b$ such that $\alpha_0$ is initial in $a$ and $\beta_0$
is terminal in $b$, else exchange
the roles of $a, b$.

Assume boundary components adjacent in $C_n$ are
separated by $c \cup x$, 
else there exists a $C_n$-allowed arc disjoint from $x, c$
and $d(x,c) \leq 2$.  For $i \notin \{1, n\}$, let
$D_i \supset w_i$ denote the closure of
the complementary component of $c \cup x$ containing $w_i$, and
let $D_1, D_{n}$ denote the closure of the union of components intersecting
$w_1$, resp.\ $w_{n}$.   Note that the $D_i$ need not be distinct, although $\mathring D_i, \mathring D_{i+1}$ are disjoint by assumption.  We show
the following claim, from which the result follows:
\begin{claim*}
  If
  $D_i \cap \alpha_0 \neq \varnothing $ and $D_{i + 1} \cap \alpha_0
  \neq \varnothing$ then there exists $y \in \mathscr U_{e_0}(a,
  c)$ such that $d(y,x) \leq 4$, and likewise for $\beta_0, b$.
\end{claim*}
\begin{proof}[Proof of claim.]
  We consider the case for $\alpha_0$; the other case is analogous.
  Fix simple disjoint 	arcs $\rho, \delta$ from
  $w_i$ and $w_{i+1}$ respectively to $\alpha_0$, both disjoint from
  $c$ and disjoint from 
  $x$ except at one endpoint; if $i = 1$, then we
  allow $\rho$ to be the point $\alpha_0 \cap a_-$.  Let $\eta'$ be the
  $C_n$-allowed concatenation of $\rho$, $\delta$, and the subarc of $\alpha_0$
  between $\rho$ and
  $\delta$.  Let $\gamma_0\subset c$ denote the subarc between
  $w_n$ and the first intersection with $\alpha_0$, or $c$ if no
  such intersection exists, and let $y$ denote the $e_0$-unicorn formed
  by concatenating an initial arc of $\alpha_0$ with $\gamma_0$, or
  $c$ if $\gamma_0 = c$.  We observe that $y, \eta'$ and $\eta', x$ both 
  intersect at most once.  Applying Proposition~\ref{prop:dist_i},
  $d(x,y) \leq d(x,\eta') + d(\eta',y) \leq 4$.
\end{proof}
\noindent
In particular, we note that $\alpha_0 \cap D_1 \neq
\varnothing$ and $\beta_0 \cap D_{n} \neq \varnothing$, and in
general $\partial D_i \cap x \neq \varnothing$ else $c$ is not
simple, hence either $D_i \cap \alpha_0 \neq \varnothing$ or $D_i \cap
\beta_0 \neq \varnothing$.
Since $n$ is odd, the
hypothesis of the claim must hold for some $i$.
\end{proof}

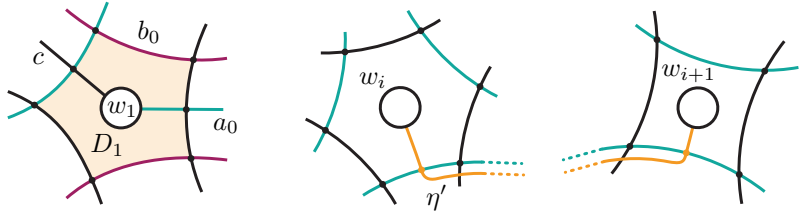
\begin{figure}[t]
\begin{tikzpicture}
\pgftext{%
\includegraphics[scale=\figscale]{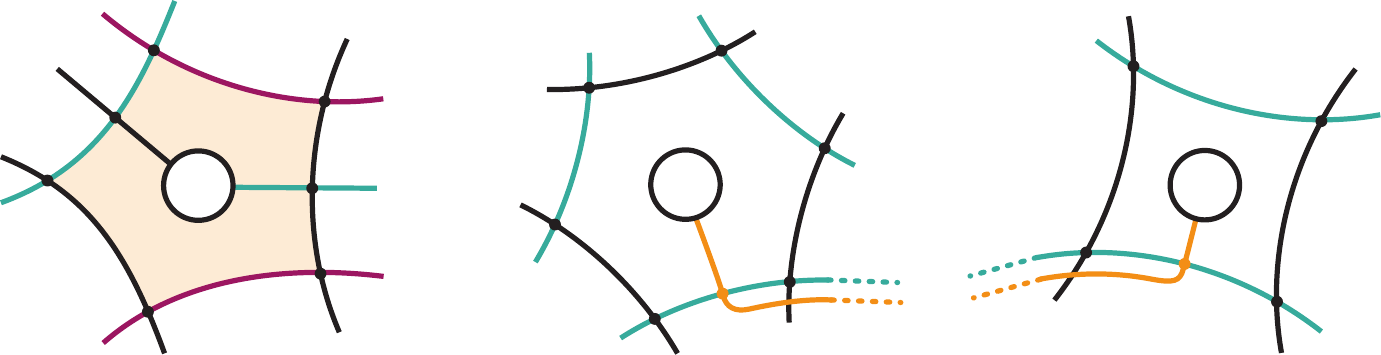}%
}%


\node at (-3.73,-.075) {$w_1$};
\node at (-3.93,-.55) {$D_1$};
\node at (-2.35,-.3) {$a_0$};
\node at (-3.375,.95) {$b_0$};
\node at (-4.85,.6) {$c$};

\node at (-.4,.3) {$w_i$};
\node at (3.8,.4) {$w_{i+1}$};

\node at (.45,-1.25) {$\eta'$};

\end{tikzpicture}
\caption{The region $D_1$; the arc $\eta'$ in the proof of the claim.}
\end{figure}

To apply Theorem~\ref{thm:gg}, it remains only to verify that
\ref{item:gg1} is satisfied:
\begin{lem}\label{lem:oddcohere}
For any disjoint $a,b \in V(\pac(\Sigma,\Gamma))$, $\diam \mathscr
U^+(a,b) \leq 7$.
\end{lem}
\begin{proof}
Let $w_i$ denote the boundary components in $C_n$, 
as in the proof of Lemma~\ref{lem:hyp_thin}.
It suffices to show that for any $x \in \mathscr U_{e_0}(\xi a,\xi b)$,
$d(x,\{a,b\}) \leq 3$.  Let $a' = \xi a, b' = \xi b$, oriented such that their
initial (resp.\ terminal) points lie in $w_1$ (resp.\ $w_n$).
We assume $x \neq a',b'$, else $d(x,\{a,b\}) \leq 1$.   
Hence, without loss of generality, let $x$ be represented by the
concatenation of subarcs $\alpha'_0 \subset a'$ and $\beta'_0
\subset b'$ such that $\alpha'_0$ is initial in $a'$ and $\beta'_0$ is
terminal in $b'$, else exchange
the roles of $a, b$.

For $i \neq \{1,n\}$, let $\rho_i'$ be the shortest path from $w_i$ to
$X = a \cup x \cup b$, and let $p_i'$ denote the point of intersection
$\rho_i \cap X$.  If $p_i' \in a, b$, let $\rho_i$ be the
concatenation of $\rho_i'$ with the shortest subarc between $p_i'$ and
$x$ along $a, b$ respectively; else let $\rho_i = \rho_i'$.  Let $p_i$
denote the endpoint of $\rho_i$ along $x$, and note that if
$p_i \in \alpha_0'$, then since $a, a'$ are disjoint
$p_i \in b \cup \alpha_0'$ and $\rho_i \cap X \subset a' \cup b$.
$a, b$ are likewise disjoint, hence $\rho_i$ is disjoint from $a$, and
analogously if $p_i \in \beta_0'$ then $\rho_i$ is disjoint from $b$.
Suppose that $p_2 \in \alpha_0'$, and let $\delta$ be the
concatenation of $\rho_2$ with the subarc between $p_2$ and $c_1$
along $\alpha_0'$.
Since $\delta$ joins boundary components
$w_1, w_2$
adjacent in $C_n$, $\delta$ is $\Gamma$-allowed; moreover, $\delta$ is
disjoint up to isotopy from $a$ and $x$, hence
$d(x,a) = 2$.  Thus assume $p_2 \notin \alpha_0'$ and, by a similar argument,
$p_{n-1} \notin \beta_0'$.

We claim that if $p_i, p_{i+1} \in \alpha_0'$, then there exists a
$\Gamma$-allowed arc $\delta$ such that $\delta$ is disjoint from $a$
and intersects $x$ at most once, and likewise for $\beta_0'$ and $b$.
In particular, in the case for $\alpha_0'$ let $\delta$ be the
concatenation of $\rho_i, \rho_{i+1}$, and the subarc along
$\alpha_0'$ between $p_i, p_{i+1}$; $\delta$ satisfies the claim, and
the remaining case for $\beta_0'$ is analogous.  Finally, since $n$ is
odd and $p_2 \in \beta_0'$ and $p_{n-1} \in \alpha_0'$, the hypotheses
of the claim are satisfied for some $i$.  Applying
Proposition~\ref{prop:dist_i},
$d(x,\{a,b\}) \leq d(x,[\delta]) + d([\delta],\{a,b\}) \leq 3$.
\end{proof}
Proposition~\ref{prop:odd_hyp} follows from Theorem~\ref{thm:gg}.

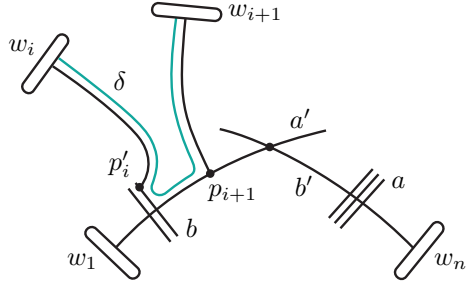
\begin{figure}[t]
\begin{tikzpicture}
\pgftext{%
\includegraphics[scale=\figscale]{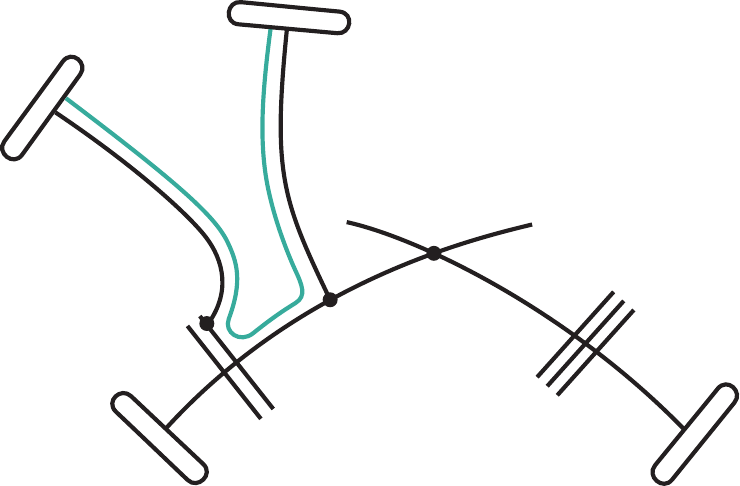}%
}%

\node at (.9,.3) {$a'$};
\node at (.95,-.6) {$b'$};
\node at (2.2,-.55) {$a$};
\node at (-.55,-1.15) {$b$};

\node at (2.9,-1.6) {$w_n$};
\node at (-2,-1.6) {$w_1$};
\node at (-2.8,1.25) {$w_i$};
\node at (.3,1.7) {$w_{i+1}$};

\node at (-1.5,.8) {$\delta$};
\node at (-1.5,-.3) {$p_i'$};
\node at (0,-.7) {$p_{i+1}$};

\end{tikzpicture}
\caption{The arc $\delta$ in the proof of Lemma~\ref{lem:oddcohere} when $p_i, p_{i+1} \in \alpha_0'$.}
\end{figure}

\section{$(\Sigma, \Gamma)$ for which $\pac(\Sigma,\Gamma)$ is not  hyperbolic}\label{sec:nonhyp}

We aim to show the converse of Theorem~\ref{thm:hyp}, namely that,
outside of some low-complexity cases, if $\Gamma$ is bipartite then
$\pac(\Sigma,\Gamma)$ is not hyperbolic. 
We first introduce the
notion of a \textit{witness subsurface}:

\begin{definition}
  An essential connected proper  
subsurface $W \subset \Sigma$ is a
  ($\Gamma$-)\textit{witness} for $\pac(\Sigma, \Gamma)$ if every
  $\Gamma$-allowed arc intersects $W$.
\end{definition}

\noindent
In the case that $\Gamma$ is bipartite we construct two independent
loxodromic actions, each supported on a disjoint witness, and thereby
deduce a quasi-isometric embedding
$\Z^2 \hookrightarrow \pac(\Sigma,\Gamma)$.  In considering such
actions, we note that the usual mapping class group does not act upon
$\pac(\Sigma,\Gamma)$ in general: if $\Gamma$ contains a non-loop edge
and is not complete, or if $\Gamma$ contains a loop but not every
loop, then there exists a mapping class whose induced permutation on
$\pi_0(\partial \Sigma) = V(\Gamma)$ does not preserve adjacency in
$\Gamma$.  We make the following definition:
\begin{definition}
  Let $\Mod(\Sigma,\Gamma) \leq \Mod(\Sigma)$ denote the subgroup of
  mapping classes $\varphi \in \Mod(\Sigma)$ 
  whose induced map $\varphi_* : \pi_0(\partial \Sigma) \to
  \pi_0(\partial \Sigma)$ defines a graph
  automorphism on $\Gamma$.
\end{definition}
Since $\Mod(\Sigma,\Gamma)$ maps $\Gamma$-allowed arcs to
$\Gamma$-allowed arcs and preserves disjointness, the following
statement is immediate:
\begin{prop}
  $\Mod(\Sigma,\Gamma)$ acts simplicially on $\pac(\Sigma,\Gamma)$. \qed
\end{prop}
\begin{rmk}
  The {pure mapping class group} $\PMod(\Sigma)$ is a subgroup of
  $\Mod(\Sigma,\Gamma)$.
\end{rmk}

\subsection{Subsurface projection}

For a witness subsurface $W \subset \Sigma$, we may extend 
(or more precisely, restrict)
the usual definitions of subsurface projection for arc graphs to
the prescribed arc graph $\pac(\Sigma,\Gamma)$.  In particular:
\begin{definition}
  Let $W \subset \Sigma$ be a $\Gamma$-witness.  Then let $\sigma_W =
  \rho_W\circ \pi :
  \pac(\Sigma,\Gamma) \to \pac(W)$ denote the \textit{subsurface
    projection} to $W$, where $\pi : \pac(\Sigma,\Gamma) \to
  \pac(\Sigma)$ is the canonical simplicial map and $\rho_W :
  \pac(\Sigma) \to \pac(W)$ is the usual subsurface projection.
\end{definition}
\begin{rmk*}\label{rmk:ssp_lip}
  We regard $\rho_W$ as a choice of true function instead of a coarse function;
  we note for $a \in \pac(\Sigma)$, the choice of image $\rho_W(a) \in
  \pac(W)$ is
  canonical up to uniformly bounded diameter
  $D$. 
Thus $\rho_W$ is
  coarsely $1$-Lipschitz, and likewise since $\pi$ is
  $1$-Lipschitz, $\sigma_W$ is
  coarsely $1$-Lipschitz, \ie\ $d_{\pac(W)}(\sigma_W(a),\sigma_W(b)) \leq d_{\pac(\Sigma,\Gamma)}(a,b)
  + 2D$ for $a,b \in V(\pac(\Sigma,\Gamma))$.  
\end{rmk*}

Crucially, subsurface projections allow us to extend loxodromic
elements in $\PMod(W)$ acting on $\pac(W)$ to loxodromics in
$\PMod(\Sigma) \leq \Mod(\Sigma,\Gamma)$.  We have the following
lemma:
\begin{lem}\label{lem:loxlift}
  Let $W \subset \Sigma$ be a $\Gamma$-witness for which there exists
  a loxodromic element $\varphi \in \PMod(W)$ acting on $\pac(W)$.
  Then any extension $\tilde\varphi \in \Mod(\Sigma,\Gamma)$ acts
  loxodromically on $\pac(\Sigma,\Gamma)$.
\end{lem}
\begin{proof}
  $\sigma_W$ is coarsely natural with respect to homeomorphisms of pairs
  $(\Sigma,W) \to (\Sigma, W)$, hence the following diagram coarsely
  commutes: 
\[
  \begin{tikzcd}
\pac(\Sigma,\Gamma) \ar[r,"\tilde\varphi_*"] \dar["\sigma_W"]
 & \pac(\Sigma,\Gamma) \dar["\sigma_W"] \\
    \pac(W) \rar["\varphi_*"] & \pac(W)
  \end{tikzcd}
\]
where $\tilde\varphi_*, \varphi_*$ are the maps induced by
$\tilde\varphi, \varphi$ respectively.  In particular, since $\varphi_*$ acts
loxodromically on $\pac(W)$ and $\sigma_W$ is coarsely $1$-Lipschitz,
$\tilde\varphi_*$ acts loxodromically on $\pac(\Sigma,\Gamma)$.
\end{proof}
\noindent
Thus we are principally interested in witnesses $W \subset \Sigma$
with pseudo-Anosov (without loss of generality, pure) mapping classes,
which act loxodromically on $\pac(W)$.  In particular, we will
typically assume $W \neq \Sigma_0^3$.

As a first application of the above, we prove that for $\Gamma' \subsetneq
\Gamma$, the canonical coarse surjection $\pi
: \pac(\Sigma,\Gamma') \to \pac(\Sigma,\Gamma)$ is almost never a
quasi-isometry.
Let $e \in E(\Gamma) \setminus E(\Gamma')$.  If $e$ is not a loop, or
$\Sigma$ has genus $g \geq 1$ or at least one
boundary component, or there exists a vertex not adjacent to the vertex
adjacent to $e$, then there exists a $\Gamma'$-witness subsurface $W \subset
\Sigma$ that is not a witness for $\pac(\Sigma,\Gamma)$ and for which
$\chi(W) \leq \chi(\Sigma) + 1$.  We observe the following:
\begin{prop}
  Let $\Gamma' \subsetneq \Gamma$ and assume that
  $\pac(\Sigma,\Gamma)$ is connected.
  If there exists a subsurface $W \neq \Sigma_0^3$ that is a
  $\Gamma'$-witness but not a $\Gamma$-witness,
  then $\pi : \pac(\Sigma,\Gamma') \to \pac(\Sigma,\Gamma)$ is not a
  quasi-isometric embedding.
\end{prop}
\begin{proof}
  Fix a $\Gamma$-allowed arc $a$ disjoint from $W$.  Let $\varphi$ be
  a pseudo-Anosov mapping class in $\PMod(W)$, and let $\tilde
  \varphi$ be an extension of $\varphi$ to $\Sigma$ by identity on
  $\Sigma \setminus W$.   Since $\varphi_*$ acts
loxodromically on $\pac(W)$, by Lemma~\ref{lem:loxlift}
$\tilde\varphi_*$ acts loxodromically on $\pac(\Sigma,\Gamma')$.  Thus
for any $b \in
V(\pac(\Sigma,\Gamma'))$, 
$\{\tilde\varphi_*^k(b)\}$ has infinite diameter in
$\pac(\Sigma,\Gamma')$, hence if $\pi$ is a quasi-isometric embedding
then likewise $\{\pi\tilde\varphi_*^k(b)\} = \{\tilde\varphi_*^k(\pi
b)\}$ has infinite diameter in $\pac(\Sigma,\Gamma)$, where the 
equality follows from the naturality of $\pi$.  But $a$ is fixed by
$\tilde\varphi$, hence $\diam \{\tilde\varphi_*^k(\pi
b)\} = 2d(a,\pi b) < \infty$, a contradiction.
\end{proof}
\noindent
In the cases listed above, if
$\pac(\Sigma,\Gamma)$ is connected and $\chi(\Sigma) \leq -3$ then the hypotheses of the proposition are satisfied:
$\chi(W) \leq -2$ and hence $W \neq \Sigma_0^3$.  
$\pi$ is not a quasi-isometry. \sqed

\subsection{Disjoint witness subsurfaces}

We show that outside of some low-complexity cases, bipartite $\Gamma$ is
equivalent to the existence of
distinct, disjoint
 witnesses $W_i$ that
support loxodromic elements of $\Mod(\Sigma,\Gamma)$.

\begin{lem}\label{lem:disjoint}
  Suppose that $\chi(\Sigma) \leq -3$ and $\Gamma$ contains an edge, and if $\Sigma = \Sigma_0^{n+1}$ then $\Gamma$ is not a $n$-pointed star.
  Then $\Gamma$ is bipartite
  if and only if there exist two disjoint, distinct
  $\Gamma$-witnesses $W_i \neq \Sigma_0^3$.
\end{lem}
\begin{proof}
  We first prove the reverse direction.  Suppose that there exist two
  disjoint $\Gamma$-witnesses, $W_1, W_2 \subset \Sigma$.  Let
  $\mathscr C_{i,j} \subset \pi_0(\partial \Sigma)$ denote the
  boundary components of $\Sigma$ contained in the $j$-th
  complementary component $C_{i,j}$ of $W_i$.  Since $W_i$ is
  essential, each $\mathscr C_{i,j}$ must be a loop-free, independent
  set in $\Gamma$, else there exists a $\Gamma$-allowed arc in
  $C_{i,j}$ disjoint from $W_i$.  Let
  $\mathscr W_1 = \pi_0(\partial W_1) \cap \pi_0(\partial \Sigma)$ and
  observe that
  $V(\Gamma) = \mathscr W_1 \sqcup \bigsqcup_j \mathscr C_{1,j}$.
  Since $W_1, W_2$ are disjoint and connected, each lies within a
  unique complementary component of the other: without loss of
  generality, assume $W_1 \subset C_{2,1}$ and $W_2 \subset C_{1,1}$,
  and observe that $W_1 \cup \bigcup_{j\neq 1} C_{1,j}$ is connected
  and disjoint from $W_2 \subset C_{1,1}$, hence likewise lies in
  $C_{2,1} \supset W_1$.  Then
  $\mathscr W_1, \bigcup_{j \neq 1} \mathscr C_{1,j} \subset \mathscr
  C_{1,1}$, hence
  $\mathscr D = \mathscr W_1 \cup \bigcup_{j \neq 1} \mathscr C_{1,j}$
  is independent in $\Gamma$ and
  $\mathscr C_{ 1,1} \sqcup \mathscr D = V(\Gamma)$ partitions
  $V(\Gamma)$ into two independent sets.

  Conversely, suppose that $V(\Gamma)$ may be partitioned into two
  independent sets $X_1, X_2$, and without loss of generality assume
  both are non-empty, else $\Gamma$ does not contain an edge.
  Moreover, we may assume $|X_1| = 1$ only if $|X_2| = 1$, and 
  $|X_2| = 1$ only if $\Gamma$ is a star with center in $X_2$, else add any isolated vertices to $X_2$.  
Let $\zeta \subset \Sigma$
  be a simple closed curve separating $X_1, X_2$ with complementary
  components $Z_1 \supset X_1, Z_2 \supset X_2$.  We observe that
  $\chi(Z_1) + \chi(Z_2) = \chi(\Sigma)$ and that $\zeta$ is essential
  if and only if neither $Z_1, Z_2$ are annuli, or equivalently, both
  $\chi(Z_1), \chi(Z_2) < 0$.  It suffices to find $\zeta$ essential
  such that some $\chi(Z_i) \leq -2$: in particular, we may choose
  $W_1 = \overline Z_i$ and $W_2$ a closed regular neighborhood of
  $\zeta$; since $\chi(W_1) \leq -2$, $W_1 \neq \Sigma_0^3$.

  If $\zeta$ is essential, then we conclude: since
  $\chi(Z_1) + \chi(Z_2) \leq -3$, at least one of
  $\chi(Z_1),\chi(Z_2) \leq -2$.  Suppose not, hence $\zeta$ is peripheral and (exactly) one of
  $Z_1, Z_2$ is an annulus.  
By our assumptions on $X_1,X_2$, in
  either case $|X_2| = 1$ and $\Gamma$ is a star, hence $\Sigma$ has
  genus.  
Fix a subsurface $Z'_2 \cong \Sigma_1^2$ containing $X_2$ and one genus.  Let $\zeta' = \partial Z'_2 \setminus \partial \Sigma$ and $Z'_1 = \overline{(\Sigma \setminus Z'_2)}$.  
Then $\zeta'$ is a simple closed curve separating $X_1, X_2$ and $\chi(Z'_1) = \chi(\Sigma) - \chi(Z'_2) = \chi(\Sigma) + 2 \leq -1$, hence $\zeta'$ is essential.
\end{proof}

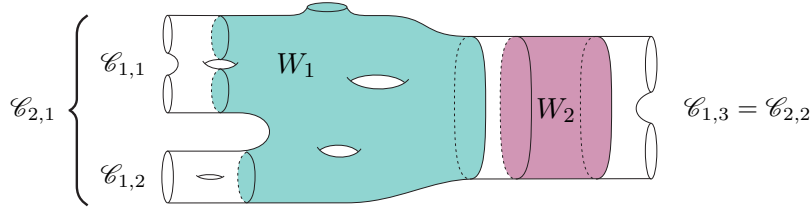
\begin{figure}[t]
\begin{tikzpicture}
\pgftext{%
\includegraphics[scale=\figscale]{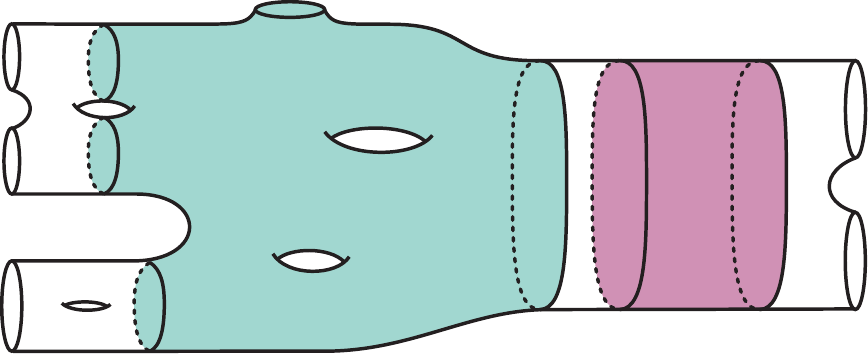}%
}%

\node at (-1.5,.5) {\large $W_1$};
\node at (1.95,-.1) {\large $W_2$};

\node at (4.5,-.1) {$\mathscr C_{1,3} = \mathscr C_{2,2}$};
\node at (-3.8,.5) {$\mathscr C_{1,1}$};
\node at (-3.8,-1) {$\mathscr C_{1,2}$};

\pgftext [at={\pgfpoint{-4.4cm}{-.1cm}}]{%
\includegraphics[scale=\figscale]{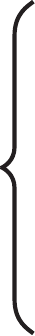}%
}%

\node at (-5,-.1) {$\mathscr C_{2,1}$};

\end{tikzpicture}
\caption{Disjoint witnesses and the collections of boundary components $\mathscr C_{i,j}$.}
\end{figure}

We conclude by showing that the existence of loxodromics supported on
disjoint witnesses implies a quasi-isometric embedding of $\Z^2$.\footnote{This statement is analogous to Exercise 3.13 in \cite{curvenotes}, mildly adapted to our setting.}

\begin{prop}\label{prop:z2embed}
  Suppose there exist distict, disjoint $\Gamma$-witnesses
  $W_1, W_2 \subset \Sigma$ that support mapping classes
  $\varphi_i \in \PMod(W_i)$ acting loxodromically on $\pac(W_i)$.
  Then $\pac(\Sigma,\Gamma)$ is not $\delta$-hyperbolic and in
  particular there exists a quasi-isometric embedding
  $\psi : \Z^2 \to \pac(\Sigma,\Gamma)$.
\end{prop}
\begin{proof}
  Let $\tau_i$ denote the translation length of $\varphi_{i*}$, and
  let $\tilde \varphi_i \in \PMod(\Sigma)$ denote the mapping class
  obtained from $\varphi_i$ by extending by identity on
  $\Sigma\setminus W_i$. Define
  $\eta : (n,m) \mapsto \tilde\varphi_{1*}^n \tilde\varphi_{2*}^m \in
  \Isom(\pac(\Sigma,\Gamma))$. Since $W_i$ is disjoint from the
  support of $\tilde \varphi_{3-i}$, $\tilde\varphi_1,\tilde\varphi_2$
  commute and $\eta$ is an action of $\Z^2$ on $\pac(\Sigma,\Gamma)$
  by isometries.  Fix a $\Gamma$-allowed arc $\alpha \subset \Sigma$
  and let $\ell_i = d(\alpha,\tilde\varphi_{i*}(\alpha))$; let
  $\psi(n,m) = \eta(n,m)(\alpha)$ denote the orbit map at $\alpha$.
  Then $d(\psi(n,m),\psi(n',m')) \leq \ell_1|n-n'| + \ell_2|m-m'|$.

  For the lower quasi-isometry bound, observe that by the disjointness
  of $W_i$ and $\supp(\tilde\varphi_{3-i})$,
  $\sigma_{W_i}(\psi(n_1,n_2))$ is coarsely equal to
  $\varphi_{i*}^{n_i}(\sigma_{W_i}(\alpha))$: in particular, there
  exists a constant $D$ such that
  $d(\sigma_{W_i} (\psi(n_1,n_2)),
  \varphi_{i*}^{n_i}(\sigma_{W_i}(\alpha))) \leq D$.  Then
  $d(\sigma_{W_i}(\psi(n_1,n_2)),\sigma_{W_i}(\psi(n_1',n_2'))) \geq
  \tau_i|n_i-n_i'| - 2D - M$, where $M > 1$ is such that
  $d(\sigma_{W_i}(\alpha),\varphi_{i*}^k(\sigma_{W_i})) \geq k\tau_i -
  M$.  Since the subsurface projections $\sigma_{W_i}$ are coarsely
  Lipschitz, the claim follows.
\end{proof}

The main theorem of this section follows immediately from
Lemma~\ref{lem:disjoint} and Proposition~\ref{prop:z2embed}:
\begin{thm}\label{thm:nonhyp}
  Suppose that $\chi(\Sigma) \leq -3$ and $\Gamma$ contains an edge,
  and if $\Sigma = \Sigma_0^{n+1}$ then $\Gamma$ is not a $n$-pointed
  star.  If $\Gamma$ is bipartite, then $\pac(\Sigma,\Gamma)$ is not
  $\delta$-hyperbolic. \qed
\end{thm}

\section{Sporadic cases}\label{sec:sporadic}

We first observe that for $\chi(\Sigma) \geq 0$ or
$E(\Gamma) = \varnothing$, $\pac(\Sigma,\Gamma)$ is either empty or a
singleton.  Similarly, we observe that $\pac(\Sigma_0^3)$ is finite
and that $\pac(\Sigma_0^3,\Gamma)$ is a full subgraph of
$\pac(\Sigma_0^3)$. Assuming $\chi(\Sigma) \leq -1$,
$E(\Gamma) \neq \varnothing$, and $\Sigma \neq \Sigma_0^3$, we
enumerate the low complexity cases not addressed in the preceding sections:

\begin{enumerate}[label=\textit{(\roman*)}]
\item \textit{Infinite diameter.} Whenever $\chi(\Sigma) \leq -2$,
  there exists a pseudo-Anosov class in $\PMod(\Sigma)$ acting
  loxodromically on $\pac(\Sigma)$, hence $\pac(\Sigma)$ is infinite
  diameter and by Theorem~\ref{thm:diam} likewise is
  $\pac(\Sigma,\Gamma)$.  Hence we need only consider
  $\Sigma = \Sigma_1^1$ with $\Gamma$ a loop, in which case
  $\pac(\Sigma_1^1,\Gamma) = \pac(\Sigma_1^1) \simeq \mathscr
  C(\Sigma_1)$ is the usual curve graph of the torus.  In particular,
  $\pac(\Sigma_1^1)$ is the Farey graph, which has infinite
  diameter. \sqed

\item \textit{Connectivity and $\delta$-hyperbolicity.}
  Theorem~\ref{thm:connected} implies the connectivity of
  $\pac(\Sigma,\Gamma)$ except when $\Sigma = \Sigma_1^2$ and $\Gamma$
  consists of only loops; if $\Gamma$ is a single loop, then
  connectivity follows from Proposition~\ref{prop:unicorn}, hence it
  remains only to consider when $\Gamma$ is two disjoint loops.
  Similarly, if $\Gamma$ is not bipartite then Theorem~\ref{thm:hyp}
  and Corollory~\ref{cor:loophyp} imply $\delta$-hyperbolicity for
  $\pac(\Sigma,\Gamma)$ except when $\Sigma = \Sigma_1^2$ and $\Gamma$
  is two disjoint loops.  We show connectedness and
  $\delta$-hyperbolicity in this case in Lemma~\ref{lem:twoloops}.

\item \textit{Non-hyperbolicity.}  If $\Gamma$ is bipartite, then
  Theorem~\ref{thm:nonhyp} implies that $\pac(\Sigma,\Gamma)$ is
  non-$\delta$-hyperbolic unless $\chi(\Sigma) \geq -2$ or
  $\Sigma = \Sigma_0^{n+1}$ and $\Gamma$ is a $n$-pointed star.  We
  show $\delta$-hyperbolicity in the latter case in
  Section~\ref{sec:star}.  For the former, it remains to consider only
  if $\Sigma = \Sigma_0^4$ or if $\Sigma = \Sigma_1^2$ and $\Gamma$ is
  a non-loop edge; in the second case, we show $\delta$-hyperbolicity
  in Lemma~\ref{lem:edge}.

  If $\Sigma = \Sigma_0^4$, then $\pac(\Sigma, \Gamma)$ is a full
  subgraph of $\pac(\Sigma_0^4)$, which is quasi-isometric to the
  Farey graph hence a quasi-tree; $\pac(\Sigma,\Gamma)$ is likewise a
  quasi-tree and thus $\delta$-hyperbolic. \sqed
 \end{enumerate}

 Collecting the results for these sporadic cases with the general
 results in Sections~\ref{sec:conn}, \ref{sec:hyp}, and
 \ref{sec:nonhyp}, we have shown Theorem~\ref{thm:conninf}:

\begingroup
\def\thethm{\ref{thm:conninf}}
\begin{thm}
  Assume that $\chi(\Sigma) \leq -1, E(\Gamma) \neq \varnothing$, and
  $\Sigma \neq \Sigma_0^3$.  Then $\pac(\Sigma, \Gamma)$ is connected
  and has infinite diameter. \qed
\end{thm}
\addtocounter{thm}{-1}
\endgroup

\noindent
We note that if $\Sigma = \Sigma^{n+1}_0$ and $\Gamma$ is a
$n$-pointed star with center $c$, then every witness subsurface must
contain $c$: $\Sigma$ does not admit disjoint $\Gamma$-witnesses.
Similarly, if $\Sigma = \Sigma_0^4$ or $\Sigma_1^2$, then $\Sigma$
does not admit two disjoint, distinct $\Gamma$-witnesses that are not
$\Sigma_0^3$.  Hence we conclude the following, proving Theorems~\ref{thm:hypclass0} and \ref{thm:hypclass}:

\begin{thm}
  Assume that $\chi(\Sigma) \leq -1, E(\Gamma) \neq \varnothing$, and
  $\Sigma \neq \Sigma_0^3$.  Then
\begin{enumerate}[label=(\roman*)]
\item \label{item:witnesseq} $\pac(\Sigma,\Gamma)$ is
  $\delta$-hyperbolic if and only if $\Sigma$ does not admit two
  distinct, disjoint $\Gamma$-witnesses that are not homeomorphic to
  $\Sigma_0^3$.
\end{enumerate}
In particular, if $\Sigma = \Sigma_0^{n+1}$ and $\Gamma$ is a
$n$-pointed star, or if $\Sigma = \Sigma_1^2$ and $\Gamma$ is a
non-loop edge, then $\pac(\Sigma,\Gamma)$ is $\delta$-hyperbolic.
Outside of these sporadic cases, \ref{item:witnesseq} is equivalent to
the following:
\begin{enumerate}[resume,label=(\roman*)]
\item $\pac(\Sigma,\Gamma)$ is $\delta$-hyperbolic if and only if
  $\Gamma$ is not bipartite. \qed
\end{enumerate}
\end{thm}

We conclude by addressing the cases when $\Sigma = \Sigma_1^2$ and
$\Gamma$ is either two loops or a non-loop edge, and when
$\Sigma = \Sigma_0^{n+1}$ and $\Gamma$ is a $n$-pointed star.

\subsection{$\Sigma = \Sigma_1^2$ and $\Gamma$ is two loops or a
  non-loop edge}

\begin{lem} \label{lem:twoloops} Suppose that $\Sigma = \Sigma_1^2$
  and $\Gamma$ is two disjoint loops.  Then $\pac(\Sigma,\Gamma)$ is
  connected and $\delta$-hyperbolic.
\end{lem}
\begin{proof}
  Let $\ell_1,\ell_2$ denote the two loops in $\Gamma$, and let $w_i$
  denote the boundary component in $\ell_i$.  Let
  $\pi_i : \pac(\Sigma, \ell_i) \to \pac(\Sigma,\Gamma)$ be the
  respective canonical map; note that $\im \pi_1, \im \pi_2$ partition
  $V(\pac(\Sigma,\Gamma))$.  We show $\pac(\Sigma,\Gamma)$ is the
  graph product of $\pac(\Sigma,\ell_1)$ and a single, non-loop edge,
  hence quasi-isometric to $\pac(\Sigma,\ell_1)$.  Since
  $\pac(\Sigma,\ell_1)$ is connected and $\delta$-hyperbolic by
  Proposition~\ref{prop:unicorn} and Corollory~\ref{cor:loophyp}
  respectively, the claim follows.

  We note that each arc $a \in \img \pi_i$ is adjacent to exactly one
  arc $a' \in \img \pi_{3-i}$.  In particular, if \eg\ $a$ is a
  $\ell_1$-allowed arc, then cutting along $a$ we obtain a three-holed
  sphere $\Sigma' = \Sigma_0^3$: there exists exactly one essential
  simple arc $a' \subset \Sigma'$ with endpoints on $w_2$, up to isotopy.
  Let $\tau : V(\pac(\Sigma,\ell_1)) \to V(\pac(\Sigma,\ell_2)) $ be
  the bijection sending each $\ell_1$-allowed arc $a$ to its unique
  disjoint $\ell_2$-allowed arc $a'$.  Since each $\pi_i$ is a graph
  inclusion, it suffices to show that $\tau$ is an isometry.  In
  particular, it is enough to find a homeomorphism of pairs
  $(\Sigma,w_1) \to (\Sigma,w_2)$ that induces $\tau$.

  For convenience, observe that we can realize $\pac(\Sigma,\Gamma)$
  by replacing the boundary components $w_i$ of $\Sigma$ with marked
  points $p_i$ on the torus.  Fix a flat metric on
  $\Sigma \cong \R^2/\Z^2$ such that $p_1 = (0,0)$ and
  $p_2 = (\theta, 0)$ for some irrational $\theta \in (0,1)$, and
  define the homeomorphism $\xi : (a,b) \mapsto (a + \theta, b)$; we
  note that $\xi : p_1 \mapsto p_2$.  For any $\ell_1$-allowed arc
  $a$, we may assume $a$ is a geodesic loop with non-zero slope and
  basepoint $p_1$, hence $\xi a$ is a geodesic loop with basepoint
  $p_2$ disjoint from $a$.  Thus $\xi a$ is $\ell_2$-allowed and by
  uniqueness $\xi a = a'$: $\xi$ induces the map $\tau$ on
  $V(\pac(\Sigma,\ell_1))$, as desired.
\end{proof}

\begin{lem} \label{lem:edge} Suppose that $\Sigma = \Sigma_1^2$ and
  $\Gamma$ is a non-loop edge.  Then $\pac(\Sigma,\Gamma)$ is
  $\delta$-hyperbolic.
\end{lem}

\begin{proof}
  As in Lemma~\ref{lem:twoloops}, without loss of generality we may
  replace $\Sigma$ with $\Sigma_{1,2}$, \ie\ replacing boundary
  components with marked points $\bar q_1, \bar q_2$ on the torus.
  Let $\Sigma' = \Sigma_{1,1}$ with marked point $q_0$; fix primitive
  generators $\alpha, \beta \in \pi_1(\Sigma',q_0)$ and a metric such
  that $\alpha, \beta$ are geodesic.  Let $p : \Sigma \to \Sigma'$ be
  the normal covering corresponding to the subgroup
  $\braket{2\alpha, \beta} < \pi_1(\Sigma',q_0)$, preserving marked
  points. We will show that $p$ defines a quasi-isometry
  $p_* : V(\pac(\Sigma,\Gamma)) \to V(\pac(\Sigma',\ell_0))$, where
  $\ell_0$ is the loop on the marked point $q_0$ and $p_*$ maps
  $\gamma \mapsto [p\hat\gamma]$, where $\hat\gamma \in \gamma$ is the
  unique (simple) geodesic representative.  Since
  $\pac(\Sigma',\ell_0) = \pac(\Sigma_1^1) \cong \mathscr C(\Sigma_1)$
  is $\delta$-hyperbolic, we conclude.

  We first verify that $p_*$ is well defined.  Let
  $\gamma \in V(\pac(\Sigma,\Gamma))$; without loss of generality,
  assume $\gamma_- = \bar q_1$.  Since $p$ is $\Pi$-injective on
  $\st(\bar q_1) \subset \Pi(\Sigma)$, $p\hat\gamma$ is essential and
  we need only check that $p\hat\gamma$ is simple.  In particular, we
  note that $V(\pac(\Sigma',\ell_0))$ corresponds bijectively with
  primitive elements in $\pi_1(\Sigma',q_0)$; let $a,b$ such that
  $[p\hat\gamma] = a\alpha + b\beta \in \pi_1(\Sigma',q_0)$, and
  assume that $p\hat\gamma$ is not simple, or equivalently, that
  $[p\hat\gamma]$ is not primitive.  Since $p\hat\gamma$ does not lift
  to a loop in $\Sigma$, $a$ must be odd, hence
  $\gcd (a,b) = k\geq 3$. Let $\omega \in \pi_1(\Sigma',q_0)$ such
  that $k\cdot \omega = [p\hat\gamma]$, hence $\hat\gamma$ is the
  concatenation of $k$ lifts of $\hat \omega$.  Since there exist
  exactly two distinct lifts of $\hat \omega$, at least two lifts
  coincide and $\hat\gamma$ is not simple, a contradiction.

  For the lower quasi-isometry bound, we observe that for any
  $\eta,\nu \in \img p_*$, the lifted arcs in $p_*^{-1}(\{\eta,\nu\})$
  are pairwise disjoint.  Therefore any path with vertices in
  $\img p_* \subset V(\pac(\Sigma',\ell_0))$ lifts to a path in
  $\pac(\Sigma,\Gamma)$ and $p_*$ is non-contracting.  For the upper
  quasi-isometry bound, it suffices to show that if $\gamma, \rho$ are
  disjoint $\Gamma$-allowed arcs, then $d(p_*\gamma, p_*\rho) \leq
  2$. We first consider the special case that $\gamma = \bar\alpha$ is
  a lift of $\alpha$ with $\bar\alpha_- = \bar q_1$. Let $\psi$ be the
  deck transformation permuting $\bar q_1,\bar q_2$.  If
  $\alpha,p_*\rho$ are not disjoint, then $\bar\alpha,\rho$ are
  disjoint but $\bar\alpha, \psi\rho$ are not. Thus
  $p_*\rho = \alpha + m\beta$ for some $m \neq 0$, hence
  $p_*\bar \alpha = \alpha$ and $p_*\rho$ are both disjoint from
  $\beta$: the claim is shown.

  We consider the induced action of $\Mod(\Sigma';q_0)\cong \SL_2(\Z)$
  on $\pi_1(\Sigma',q_0)$; let
  $H = \braket{2\alpha,\beta} < \pi_1(\Sigma',q_0)$, and let
  $\stab(H) < \Mod(\Sigma';q_0) $ denote the subgroup of mapping
  classes $\varphi$ such that $\varphi_*(H) =
  H$. $V(\pac(\Sigma',\ell_0))$ corresponds to the disjoint union of
  the sets of primitive elements in $H$ and
  $H^c \dfn \pi_1(\Sigma_{1,1},q_0) \setminus H$ respectively, and we
  note that $\stab(H)$ acts by isometries on both sets.  Moreover,
  $\img p_*$ is precisely the set of primitive elements in $H^c$.  To
  show that $p_*$ is $2$-Lipschitz in general, it suffices to show
  that $\stab(H)$ acts transitively on $\img p_*$.  In particular, if
  $\gamma$ is a $\Gamma$-allowed arc in $\Sigma$, then let
  $\varphi \in \stab(H)$ such that
  $\varphi: p_*\gamma \mapsto \alpha$. Then there exists a lift
  $\bar \varphi : \Sigma \to \Sigma$ such that
  $\bar \varphi : \hat \gamma \to \bar \alpha$; noting that
  $\bar\varphi$ acts isometrically on $\pac(\Sigma,\Gamma)$ and
  $p_*$-intertwines the isometric action of $\varphi$ on
  $\pac(\Sigma',\ell_0)$, we reduce to the special case
  above. Finally, we verify transitivity.  Let
  $\omega = a\alpha + b\beta \in H^c$; we note that $a,b$ are coprime
  and $a$ is odd.  Since $a,b$ are coprime, there exist $c,d \in \Z$
  such that $ad - bc = 1$.  Moreover, we may choose $c$ to be even,
  else replace $c$ with $c + a$ and $d$ with $d + b$.  Let
  \[\varphi = \begin{pmatrix} a & c\\ b & d
    \end{pmatrix} \in \Mod(\Sigma';q_0) \cong \SL_2(\Z),\] noting that
  $\det(\varphi) = ad - bc = 1$ and that since $a$ is odd and $c$
  even, $\varphi \in \stab(H)$.
  $\varphi(\alpha) = a\alpha + b\beta = \omega$, hence maps
  $\alpha \mapsto \omega$ as desired.

  It remains only to check that $p_*$ is coarsely surjective.  We show
  that $\stab(H)$ acts transitively on primitives in $H$, hence
  $\stab(H)\setminus V(\pac(\Sigma,\ell_0))$ has exactly two vertices,
  one of which has fiber $\img p_*$: $\img p_*$ is $2$-coarsely dense
  in $V(\pac(\Sigma,\ell_0))$.  Let $\nu = c\alpha + d\beta \in H$ be
  primitive, hence $c$ is even and there exist $a, b$ such that
  $ad - bc = 1$.  Thus $ad$, and likewise $a$, must be odd.  Let
\[\xi = 
\begin{pmatrix} a & c \\ b & d
\end{pmatrix}.
\]
As above, we may verify that $\xi \in \stab(H)$ and $\xi : \beta \mapsto \nu$.
\end{proof}

\subsection{$\Sigma = \Sigma_0^{n+1}$ and $\Gamma$ is a $n$-pointed
  star, $n \geq 3$.}\label{sec:star}

The case when $\Sigma$ is a $(n+1)$-holed sphere and $\Gamma$ is a
$n$-pointed star is analogous to the ray graph in \cite{bavard}, but
rather than a Cantor set of punctures we consider only finitely many.
We reproduce Bavard's argument for $\delta$-hyperbolicity below,
suitably simplified for our purposes.

\begin{lem}
  Let $\Sigma = \Sigma_0^{n+1}$ and $\Gamma$ a $n$-pointed star with
  center $c$.  Then $\pac(\Sigma,\Gamma)$ is quasi-isometric to
  $\pac(\Sigma,\ell_0)$, where $\ell_0$ is a loop on $c$.
\end{lem}
\begin{proof}
  It suffices to define a quasi-isometry
  $\xi : V(\pac(\Sigma,\Gamma)) \to V(\pac(\Sigma,\ell_0))$. Given an
  isotopy class of $\Gamma$-allowed arcs
  $a \in V(\pac(\Sigma,\Gamma))$, we will always assume an orientation
  such that $a_- = c$.  Let $N$ be a regular neighborhood of
  $\alpha \cup a_+$ for $\alpha \in a$ a choice of representative; let
  $\hat a$ denote the isotopy class of $\partial N$.  Since
  $n \geq 3$, $\hat a$ is essential; $\partial \hat a = \{c\}$, hence
  $\hat a$ is an isotopy class of $\ell_0$-allowed arcs.  Note that
  $\hat a$ is independent of the choice of $\alpha$ and $N$, and let
  $\xi : a \mapsto \hat a$.

If $\alpha\in a, \beta\in b$ are disjoint $\Gamma$-allowed arcs, then
we may choose regular neighborhoods of $\alpha \cup a_+$ and $\beta
\cup b_+$ respectively such that $i(\hat a, \hat b) \leq 2$.  Thus by
Proposition~\ref{prop:dist_i} $d(\xi a,\xi b) \leq 3$,
and in general $\xi$ is $3$-Lipschitz.  
For the lower quasi-isometry bound, consider $a, b \in
V(\pac(\Sigma,\Gamma))$ such that $d(\xi a,\xi b) = m$, and let $u_1 =
\xi a, u_2,\ldots,u_m=\xi b$  be distinct $\ell_0$-allowed arcs
forming a geodesic path in $\pac(\Sigma,\ell_0)$.  We construct a path
$v_1 = a, v_2, \ldots, v_m = b$ of $\Gamma$-allowed arcs in
$\pac(\Sigma,\Gamma)$, hence $d(a,b) \leq d(\xi a,\xi b)$, as desired.

We note that each $u_i$ separates $\Sigma$ into two complementary
components homeomorphic to $\Sigma_{0,1}^k, \Sigma_{0,1}^{n-k}$
respectively, and that for $i \neq 1,m$, both $u_{i-1},u_{i+1}$ must
lie in the same component: since $(u_i)$ is geodesic, $u_{i-1},
u_{i+1}$ cannot be adjacent in $\pac(\Sigma,\ell_0)$, hence must
intersect.  Let $S_i \subset \Sigma$ denote the complementary
component of $u_i$ that contains $u_{i-1}$ or $u_{i+1}$; let $S_i'$
denote the the other component, which contains neither.  
We note that $S_i$ must contain at
least two boundary components of $\Sigma$ (that are not $c$), else
every $\ell_0$-allowed arc in $S_i$ belongs to $u_i$ and $u_{i-1} =
u_i$ or $u_{i+1} = u_i$, a contradiction since each $u_i$ is distinct.
Similarly, $S_i'$ contains at least one boundary component $d_i \neq
c$, since $u_i$ is essential; we moreover claim that $S_i', S_{i+1}'$
are disjoint, else \eg\ $u_i = \partial S_i' \subset S_{i+1}'$, a
contradiction with our choice of $S_{i+1}'$.
Finally, choose $v_i$ to be an arc between $c$ and $d_i$ in $S_i'$.
Since $S_i'$ is disjoint from $S_{i+1}'$, $v_i$ is disjoint from
$v_{i+1}$, and note that we may choose $v_1 = a$ and $v_m = b$ since
the complementary component of $u_1$ containing $a$ contains only one
boundary component of $\Sigma$, hence $a \not\subset S_1$, and
likewise $b \not\subset S_m$.  $(v_i)$ is the desired path.

Coarse surjectivity follows from a similar argument: as above, a
complementary component $S$ of an (essential) $\ell_0$-allowed arc $u$
must contain at least one component $d \neq c$ of $\partial \Sigma$.
Any choice of simple arc $\alpha$ from $c$ to $d$ in (the interior of)
$S$ is $\Gamma$-allowed and has a regular neighborhood $N \subset S$.
Hence $a = [\alpha] \in V(\pac(\Sigma,\Gamma))$ and $\xi a = \hat a$
is disjoint from $u$: $\xi$ is $1$-coarsely surjective.
\end{proof}
By Corollary~\ref{cor:loophyp} $\pac(\Sigma,\ell_0)$ is
$\delta$-hyperbolic, hence we may conclude that $\pac(\Sigma,\Gamma)$
is likewise $\delta$-hyperbolic. \sqed

\section{Acknowledgements}

The author would like to thank Mladen Bestvina for numerous helpful discussions.  This work was supported in part by NSF award no.\ 1840190: \textit{RTG: Algebra, Geometry, and Topology at the University of Utah}.

\bibliographystyle{amsalpha}
\bibliography{pac}

\end{document}